\documentclass[11pt]{article}
\usepackage{fancyhdr,fancybox}

\usepackage{graphicx}
\usepackage{caption2}
\usepackage{amsfonts,amsmath}
\usepackage{natbib}

\usepackage{floatrow}
\newcommand{\un}{{\bf 1}}
\newcommand{\bR}{{\mathbb{R}}}
\newcommand{\bE}{{\mathbb{E}}}
\newcommand{\lck}{\left\lbrack}
\newcommand{\rck}{\right\rbrack}
\newcommand{\lce}{\left\lbrace}
\newcommand{\rce}{\right\rbrace}

\newcommand{\Un}{{\bf 1}}
\newcommand{\dd}{\, \mbox{d}}

\newtheorem{proposition}{Proposition}

\begin{document}
\thispagestyle{empty}

\begin{center}
{\Large
\textbf{Estimating second-order characteristics of}

 \medskip

\textbf{inhomogeneous spatio-temporal point processes:}} \\

\medskip

\textbf{{\large
influence of edge correction methods and intensity estimates}}

\bigskip
\bigskip

{\large Edith Gabriel}

\bigskip

{\small
Laboratoire de Math\'{e}matiques, Universit\'{e} d'Avignon,

33 rue Louis Pasteur, 84000 Avignon, France

\medskip

edith.gabriel@univ-avignon.fr
}

\end{center}

\bigskip

\textbf{Abstract:}
We restrict our attention to space-time point pattern data for which we have a single realisation within a finite region. Second-order characteristics are used to analyse the spatio-temporal structure of the underlying point process. In particular, the space-time inhomogeneous pair correlation function and $K$-function measure the spatio-temporal clustering/regularity and the spatio-temporal interaction, and thus help with model choice. Non-parametric estimators of the second-order characteristics require information from outside the study region, resulting to the so-called edge effects which have to be corrected, and depend on first-order characteristics, which have to be estimated in practice.
Here, we extend classical edge correction factors to the spatio-temporal setting and compare the performance of the related estimators for stationary/non-stationary and/or isotropic/anisotropic point patterns. Then, we explore the influence of estimated intensity function on these estimators.

\medskip

\textbf{Keywords}: Edge correction; Intensity estimation; Point process;
Reduced second moment measure; Spatio-temporal data.

\section{Introduction}
\label{intro}

Space-time point pattern data (also called spatio-temporal point pattern data) are increasingly available in a wide range of scientific settings \citep{verejones2009}. Data-sets of this kind often consist of a single realisation of the underlying process.  Usually, separate analyses the spatial and the temporal components are of limited value, because
the scientific objectives of the analysis are to understand and to model the underlying spatio-temporally interacting stochastic mechanisms.
Generic methods for analysing such processes are growing; see for example \cite{diggle2006}, \cite{diggle2010} and \cite{cressie2011}. There is already an extensive literature on the use of point process models in the specific field of seismology; see, for example, \cite{zhuang2002} and references therein.
There are basically two ways for modeling space-time point patterns \citep{diggle2006,verejones2009}. The first is descriptive and aims at providing an empirical description of the data, especially from second-order characteristics. The second is mechanistic and aims at constructing parametric point process models by specifying parametric models for the conditional intensity function. Here, we will consider the former and analyses will be based on extensions of the pair correlation function and the Ripley's $K$-function \citep{ripley1977} to summarize a spatio-temporal point pattern and test hypotheses about it.

In many applications the point patterns cannot be considered homogeneous.
This is particularly so in epidemiological studies where the observed point pattern is spatially and temporally inhomogeneous, as the pattern of incidence of the disease reflects both the spatial distribution of the population at risk and systematic temporal variation in risk \citep{gabriel2009}.
In the spatial case, there are different approaches to deal with the inhomogeneity  \citep{illian2008}: some consist in considering subplots of the study region and employing means of the methods for homogeneous point processes \citep{allard2001, brix2001}, others have their own underlying theory, mainly based on the inhomogeneous $K$-function proposed by \cite{baddeley2000}. The latter have been for instance used in environmental epidemiology \citep{diggle2007}, ecology \citep{law2009} or economy \citep{arbia2012}. The inhomogeneous $K$-function has been extended to the spatio-temporal setting by \cite{gabriel2009}. Second-order characteristics  are thus analysed from the spatio-temporal inhomogeneous $K$-function (STIK-function) or equivalently from the spatio-temporal pair correlation function under the assumption of second-order intensity re-weighted stationarity \citep{diggle2010,gabriel2010,gabriel2013}. Spatio-temporal separability of the STIK-function has been studied in \cite{moller2012}.

Here, we restrict our attention to space-time point processes for which we have a single realisation within a finite (spatio-temporal) region. Estimating the STIK-function or the pair correlation function thus faces two issues. On the one hand, they depend on first-order characteristics (see Section~\ref{sec:1}) which are unknown in practice, but replacing the intensity by an estimate must be made carefully as it may imply bias \citep{baddeley2000}.
The problems encountered when using the same point pattern to estimate both
a spatially varying intensity and second-order characteristics can be overcome by adjusting the estimate of intensity so as to take account of explanatory variables; see for instance \cite{diggle2007}
in the context of case-control data or \cite{mrkvicka2012} for the estimation of clustered point processes with inhomogeneous cluster centres. As suggested in \cite{arbia2012} and \cite{diggle2007}, other ways to get around this difficulty is to assume a parametric model for the intensity \citep{moller2003,liu2007} or to allow separation of first- and second-order effects. The latter is particularly difficult to check in the absence of independent replications, as we cannot in general make an empirical distinction between first-order and second-order effects without additional assumptions.
On the other hand, estimating the STIK-function or the pair correlation function requires information from outside the study region, resulting to the so-called edge-effects which have to be corrected. In the spatial case, the four main edge correction methods are: the border method (all sample points that are closer to the border than to their nearest neighbours are eliminated), the guard method (the sample points within the guard zone are not used for the analysis, but used as some nearest neighbours), the toroidal method (it makes eight identical copies of the rectangular study region and places them around the original one), the weighted method (or the isotropic method), also called Ripley's method (each weight is equal to the probability that points around the location of sample point $i$ will be in the study region). There is an extensive literature on edge correction for inter-points distance methods in the spatial case; see for instance \cite{baddeley1999}, \cite{cressie1993}, \cite{diggle2003}, \cite{illian2008}, \cite{law2009}, \cite{li2007} and \cite{ripley1988} for the definitions of usual edge correction factors, \cite{goreaud1999}, \cite{haase1995} and \cite{pommerening2006} for some other specific ones, and \cite{yamada2003} for an empirical comparative study in the homogeneous case. In the spatio-temporal case, \cite{cronie2011} propose some edge correction methods when considering a given parametric model for the point pattern. Finally, \cite{baddeley1999} points out that edge effects are more difficult to overcome in dimension greater than two.

As above mentioned, estimating first- and second-order characteristics from a single realisation of point process and correcting edge effects have been widely studied in the spatial case. The aim of this paper is
to discuss and explore the influence of these two issues in the spatio-temporal setting, particularly on non-parametric estimates of the STIK-function and the pair correlation function.
The structure of the remainder of the article is as follows.
Section~\ref{sec:1} reviews some existing results concerning the second-order
characteristics of inhomogeneous spatio-temporal point processes and their estimation. In Section~\ref{sec:2}, we extend classical edge correction factors to the spatio-temporal setting and provide an empirical comparative study of the related estimators. Section~\ref{sec:3} discusses the influence of using an estimate of the intensity on the performance of the estimators from a simulation study. Section~\ref{sec:4} is a short discussion.

\section{Description of second-order characteristics}
\label{sec:1}

\subsection{Spatio-temporal point processes}

Following \cite{diggle2010} and \cite{gabriel2009}, let us consider spatio-temporal point processes whose events are defined as countable sets of points ${x=(s, t)}$ where $s$ denotes location and $t$ denotes time. In practice we observe $n$ events $x_i=(s_i,t_i)$ within a bounded spatio-temporal region $S \times T \subset \bR^2 \times \bR$.
Here, we focus on the case where additional variables describing characteristics of the events (marks) are not available.
In the following, $N(A)$ denotes the number of events in an arbitrary region $A$. For formal definitions of point process characteristics, see for example \cite{daley2003} or \cite{moller2003}.

\subsection{First-order and second-order characteristics}
\label{subsec:intensity}

First-order characteristics are described by the intensity of the process,
$$\lambda(x)= \lambda(s,t) = \lim_{|d s| \to 0, |dt| \to 0}
  \frac{\bE \lck N(d s \times dt) \rck}{|d s| |dt|}, $$
where $d s$ defines a small spatial region around the location $s$, $|d s|$
is its area, $dt$ is a small interval containing the time $t$, $|dt|$ is its
length. If the intensity is constant, $\lambda(s,t) = \lambda$ for all $(s,t)$, then the process is called homogeneous or first-order stationary.

\noindent
The relationship between numbers of events in pairs of subregions within $S \times T$
is described by the second-order characteristics.
The second-order intensity is defined as
$$  \lambda_2 (x,x^\prime) = \lim_{|A|,
    |A^\prime| \to 0} \frac{\bE \lck N(A) N(A^\prime) \rck}
  {|A||A^\prime| }, $$
  where $A = d s \times dt$ and $A^\prime= d s^\prime
\times dt^\prime$ are small cylinders containing the points $x$
and $x^\prime$ respectively, $x \neq x^\prime$.

\noindent
Equivalent descriptors of second-order characteristics include the
covariance density,
$$ \gamma(x,x^\prime) =  \lambda_2 (x,x^\prime) - \lambda(x) \lambda(x^\prime)$$
and the radial distribution function or the
point-pair correlation function \citep{cressie1993,diggle2003}:
 \begin{equation}
 \label{eq:pcf}
 g(x,x^\prime) = \frac{\lambda_2 (x, x^\prime)} {\lambda (x) \lambda (x^\prime)}.
 \end{equation}
For a spatio-temporal Poisson process, the covariance density is identically zero and the pair correlation function is identically 1. Larger or smaller values than these benchmarks therefore indicate informally how much more or less likely it is that a pair of events will occur at the specified
locations than in a Poisson process with the same intensity.

\noindent In practice, it is often difficult to estimate these moments and one may need some relaxing assumptions.

\subsection{Spatio-temporal stationarity}

A point process is first-order and second-order stationary in both space and time, if:
$$
\lambda(s,t) = \lambda \ \text{ and } \ \lambda_2 \big((s,t),(s^\prime,t^\prime)\big) = \lambda_2(s-s^\prime,t - t^\prime).$$
A stationary spatio-temporal point process is isotropic if $\lambda_2 \big((s,t), (s^\prime,t^\prime)\big) = \lambda_2 (u,v)$, where $u= \| s - s^\prime \|$ and $v = |t - t^\prime|$ denote spatial and temporal distances respectively.

\noindent
\cite{baddeley2000} introduce a weaker assumption in the purely spatial case, called second-order intensity reweighted stationarity. \cite{gabriel2009} extend it to the spatio-temporal setting as follows: an isotropic spatio-temporal point process is second-order intensity reweighted stationary if its
its pair correlation function depends only on the spatio-temporal difference vector $(u,v)$.

\subsection{Separability}

A spatio-temporal point process is first-order separable if its intensity
$\lambda(s,t)$ can be factorised as $$\lambda(s,t) = \lambda_1(s) \lambda_2(t), \ \text{for all } (s,t) \in S \times T,$$
where $\lambda_1$ and $\lambda_2$ are non-negative functions.

\noindent
The pair correlation function is separable in space and time if
$$g\big( (s,t),(s',t') \big) =g_1(s,s')g_2(t,t'),$$
where $g_1$ and $g_2$ are non-negative functions. Then for an isotropic second-order intensity reweighted stationary  point process, we have
$$g\big( (s,t),(s',t') \big) =g_1(u)g_2(v).$$
Some diagnostic procedures for checking hypotheses of second-order spatio-temporal separability are proposed in \cite{moller2012}.

\subsection{Spatio-temporal inhomogeneous $K$-function}
\label{subsec:STIK}

Second-order characteristics can also be described through an extension of Ripley's $K$-function.
For a second-order intensity reweighted stationary, isotropic, spatio-temporal point process, \cite{gabriel2009} define the space-time inhomogeneous  $K$-function by
\begin{equation}
  \label{eq:STIKfunction1}
K(u,v) =  2 \pi \int_{-v}^v \int_0^u g(u',v') \ u' \dd u' \dd v'.
\end{equation}

For an inhomogeneous spatio-temporal Poisson process with intensity $\lambda(x)$, the second-order intensity is $\lambda_2(x,x^\prime) = \lambda(x) \lambda(x^\prime)$, hence the pair correlation function $g(u,v)$ is identically 1 and the STIK-function is $K(u,v) = 2 \pi u^2 v$.
Thus, $K(u,v) - 2 \pi u^2v$ can be used as a measure of the spatio-temporal aggregation or regularity, using an inhomogeneous Poisson process as a benchmark. Positive values of $K(u, v) - 2\pi u^2v$ indicate clustering, or aggregation, at spatial and temporal separations less than $u$ and $v$, respectively, while negative values indicates regularity.

The STIK-function can also be used as a measure of spatio-temporal interaction.
Its separability into purely spatial and temporal components indicates absence of interaction. In this case the ratio
$K(u,v)/\{K_S(u) K_T(v)\} $ is constant \citep{moller2012}. It is identically 1 for a Poisson process.

\subsection{Estimation}

An unbiased estimator of the STIK-function, proposed in \cite{gabriel2009} under the assumption of isotropy, is
\begin{equation}
  \label{eq:STIKhat1}
  \widehat {K}(u,v) = \sum_{i=1}^{n}
  \sum_{j\neq i} \frac{1}{w_{ij}} \frac{1}{\lambda(x_i) \lambda(x_j)}{\bf 1}_{\lce \|s_i - s_j\| \leq u
  \ ; \ |t_i - t_j| \leq v \rce},
\end{equation}
where $w_{ij}$ is an edge correction factor to deal with spatial-temporal edge effects. The proof of unbiasedness is given in Appendix 1 for the edge correction factors defined in Section~\ref{subsec:edge}.

The space-time pair correlation function defined in Equation~(\ref{eq:pcf}) can be estimated by
\begin{equation}
\label{eq:ghat}
\widehat g(u,v) = \frac{1}{4 \pi u} \sum_{i=1}^{n}
  \sum_{j \neq i} \frac{1}{w_{ij}} \frac{k_{s}(u-\|s_i-s_j \|)k_{t}(v-|t_i-t_j|)}{\lambda(x_i) \lambda(x_j)},
\end{equation}
where $w_{ij}$  is defined in Equation~(\ref{eq:STIKhat1}) and $k_{s}(\cdot)$, $k_{t}(\cdot)$ are kernel functions with bandwidths $h_s$ and $h_t$. Experience with pair correlation function estimation recommends box kernels, see \cite{illian2008}.

\medskip

 These estimators depend on an edge correction factor and assume that the intensity is known. In practice the intensity function is unknown and have to be estimated. So, the questions are: which method can be used to correct edge effects and can be used to estimate the intensity function? what are their influence on the performance of the STIK-function and the pair correlation function? This is explored in the following sections.

\section{Influence of edge correction methods}
\label{sec:2}

The literature on edge correction methods is more extensive in the spatial case, see e.g. \cite{baddeley1999}, \cite{illian2008}, \cite{ripley1988} and \cite{stoyan1994}, than in higher dimensions \citep{baddeley1993,cronie2011,diggle1991,mamaghani2010}.

In this section we extend three classical spatial edge correction factors to the spatio-temporal setting and compare the performance of the related estimators of the second-order characteristics for stationary/non-stationary and/or isotropic/anisotropic point patterns.

\subsection{Edge correction methods}
\label{subsec:edge}

The simplest approach to deal with spatio-temporal edge effects consists in correcting them separately \citep{diggle1995}. Thus, the edge correction factor in Equations~(\ref{eq:STIKhat1}) and (\ref{eq:ghat}) is the product of a spatial edge correction factor and a temporal edge correction factor.
Here, we consider spatial edge correction factors (as described in \cite{baddeley1999} and \cite{ripley1988}) running for any shape of $S$. This eliminates in particular the toroidal method.

In the study of Section~\ref{subsec:simu} we also consider the case where no edge correction is performed:
$w_{ij} = |S \times T|.$

\subsubsection{Isotropic edge correction method}

The weight is proportional to the product between the Ripley edge correction factor \citep{ripley1977} and its analog in one-dimension:
$$w_{ij} = |S \times T| w_{ij}^{(t)} w_{ij}^{(s)},$$
where the temporal edge correction factor $w_{ij}^{(t)} = 1$ if both ends of the interval of length $2 |t_i - t_j|$ centred at $t_i$ lie within $T$ and $w_{ij}^{(t)}=1/2$ otherwise \citep{diggle1995} and
$w_{ij}^{(s)}$ is the proportion of the circumference of a circle centred at the location $s_i$ with radius $\|s_i -s_j\|$ lying in $S$.

\subsubsection{Border and modified border methods}

These methods restricts attention to those events lying more than $u$ units away from the boundary of $S$  \citep{diggle1979} and more than $v$ units away from the boundary of $T$. Thus, for $d(s_i,S)$ denoting the distance between $s_i$ and the boundary of $S$ and $d(t_i,T)$ denoting the distance between $t_i$ and the boundary of $T$, we have for the border method
$$w_{ij} = \dfrac{\sum_{j=1}^n \Un_{\{d(s_j,S) > u \ ; \ d(t_j,T) >v\}} / \lambda(x_j)}{\Un_{\{ d(s_i,S) > u \ ; \ d(t_i,T) >v \}}}$$
and
$$w_{ij} = \dfrac{|S_{\ominus u}|\times|T_{\ominus v}|}{\Un_{\{ d(s_i,S) > u \ ; \ d(t_i,T) >v \}}}$$
for the modified border method, where $S_{\ominus u}$ and $T_{\ominus v}$ are the eroded spatial and temporal region respectively, obtained by trimming off a margin of width $u$ and $v$ from the border of the original region \citep{baddeleyturner2000}.

\subsubsection{Translation method}

The weight is the proportion of translations of $(x_i,x_j)$ which have both $x_i$ and $x_j$ inside $S \times T$ ; see \cite{osher1983} for the spatial case:
$$w_{ij} =|S \cap S_{s_i-s_j}| \times |T \cap T_{t_i-t_j}|,$$
where $S_{s_i-s_j}$ and $T_{t_i-t_j}$ are the translated spatial and temporal regions.

\subsection{Simulations}
\label{subsec:simu}

We examine the performance of the STIK-function and the pair correlation function via Monte Carlo simulations.
\cite{gabriel2013} covers many of the models encountered in applications of point process methods to the study of spatio-temporal phenomena. In particular an inhomogeneous Poisson process is defined by:
\begin{enumerate}
\item The number $N(S \times T)$ of events within the region $S \times
  T$ follows a Poisson distribution with mean $\int_S \int_T
  \lambda(s,t) \dd t \dd s$.
\item Given $N(S \times T) = n$, the $n$ events in $S \times T$ form
  an independent random sample from the distribution on $S \times T$ with
  probability density function $f(s,t) \! = \! \lambda(s,t) / \!\! \int_S \! \int_T \!
  \lambda(s', t') \dd t' \! \dd s'$.
\end{enumerate}
and a Poisson cluster process is defined by:
\begin{enumerate}
\item Parents form a Poisson process with intensity $\nu(x)$.
\item The number of offspring per parent is a random variable $N_c$ with
  mean $m_c$, realised independently for each parent.
\item The positions and times of the offspring relative to their parents are independently and identically distributed according to a trivariate
  probability density function $\varphi(\cdot)$ on $\bR^2 \times \bR$.
\item The final process is composed of the superposition of the offspring
  only.
\end{enumerate}
Realisations of these point processes can be generated using the {\tt R} \citep{R} package {\tt stpp} \citep{gabriel2013}.

\bigskip

\begin{figure}[h]
 \begin{center}
  \includegraphics[width=3.75cm,height=3.75cm]{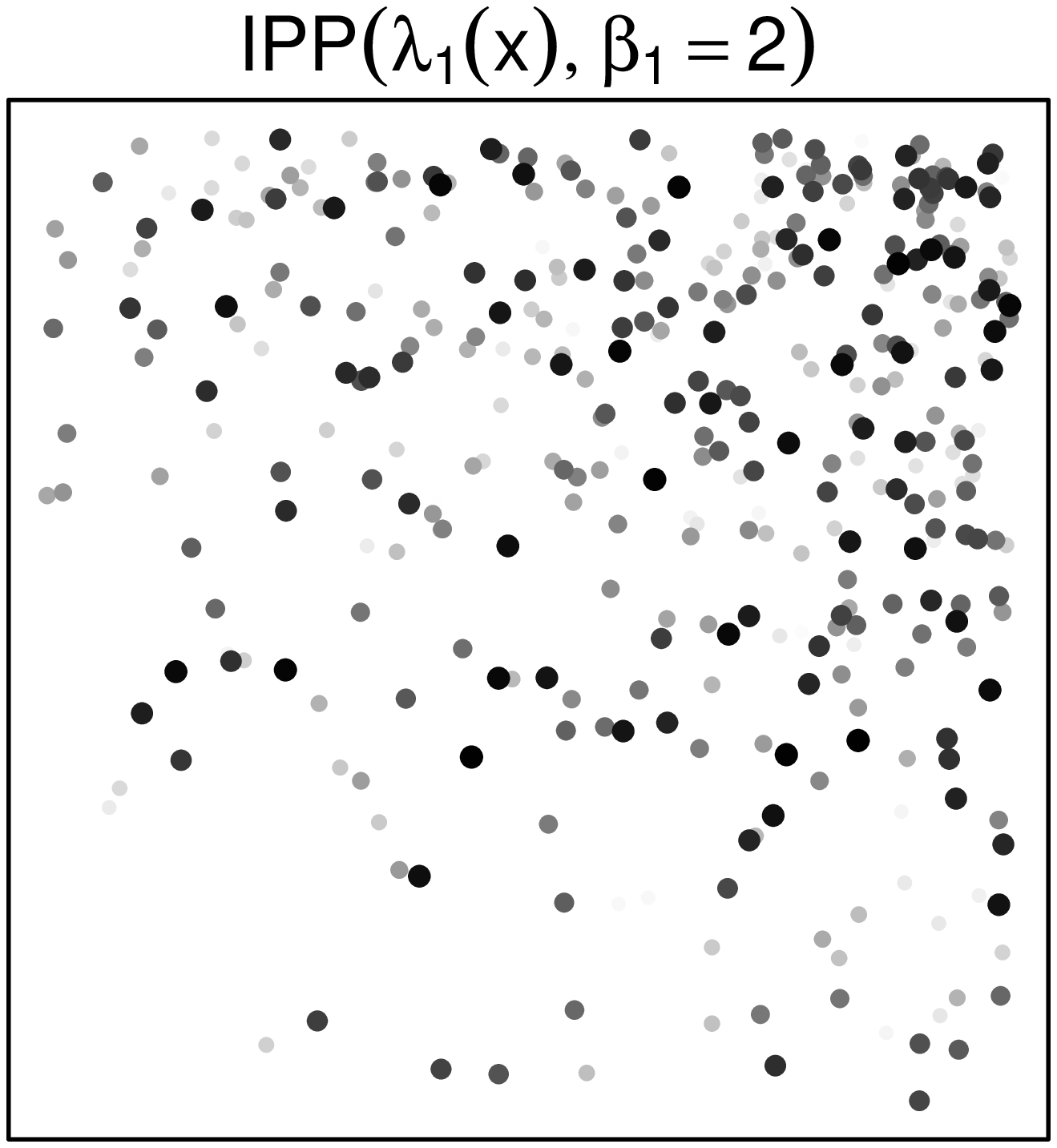}
    \includegraphics[width=3.75cm,height=3.75cm]{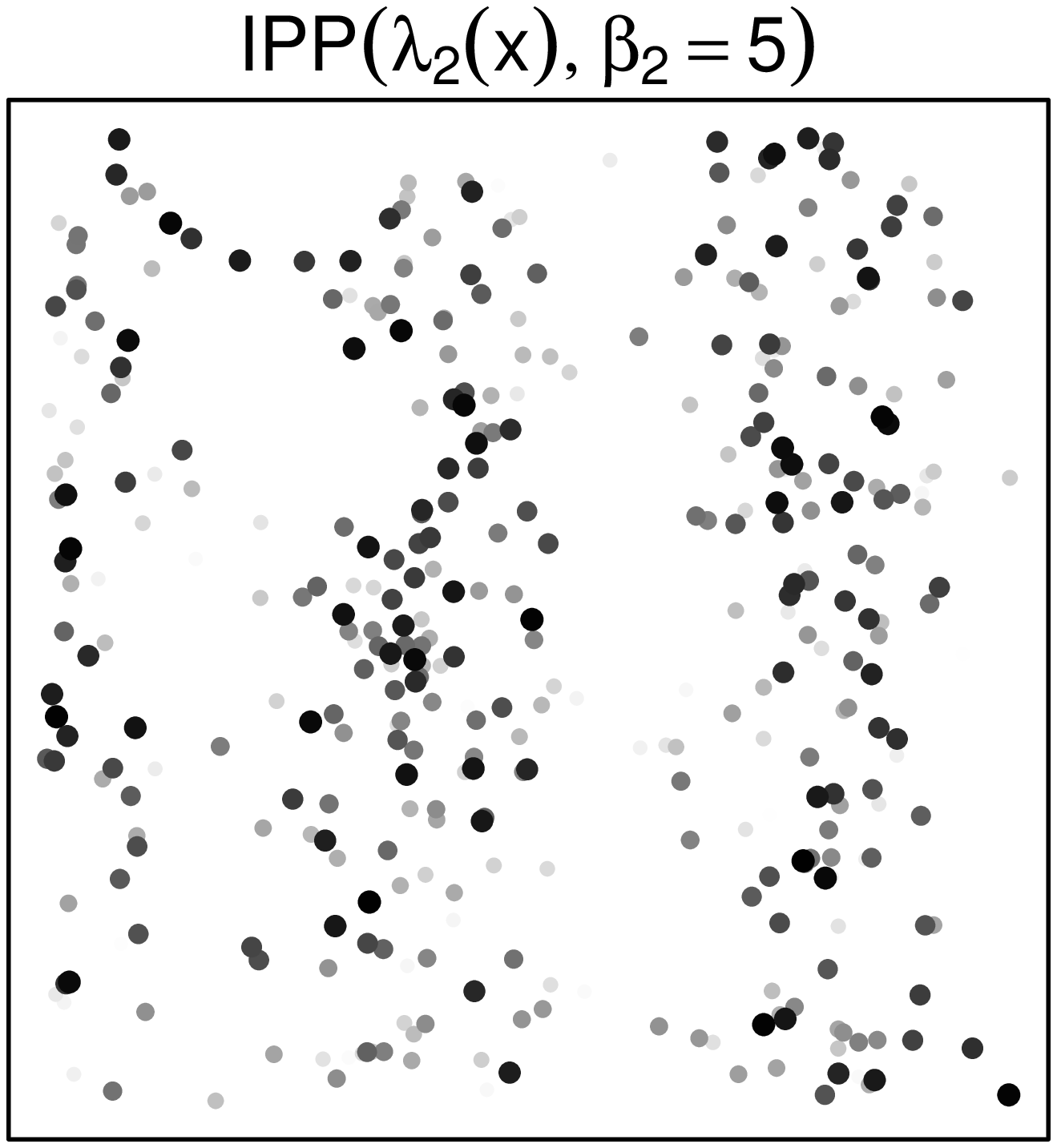}
  \includegraphics[width=3.75cm,height=3.75cm]{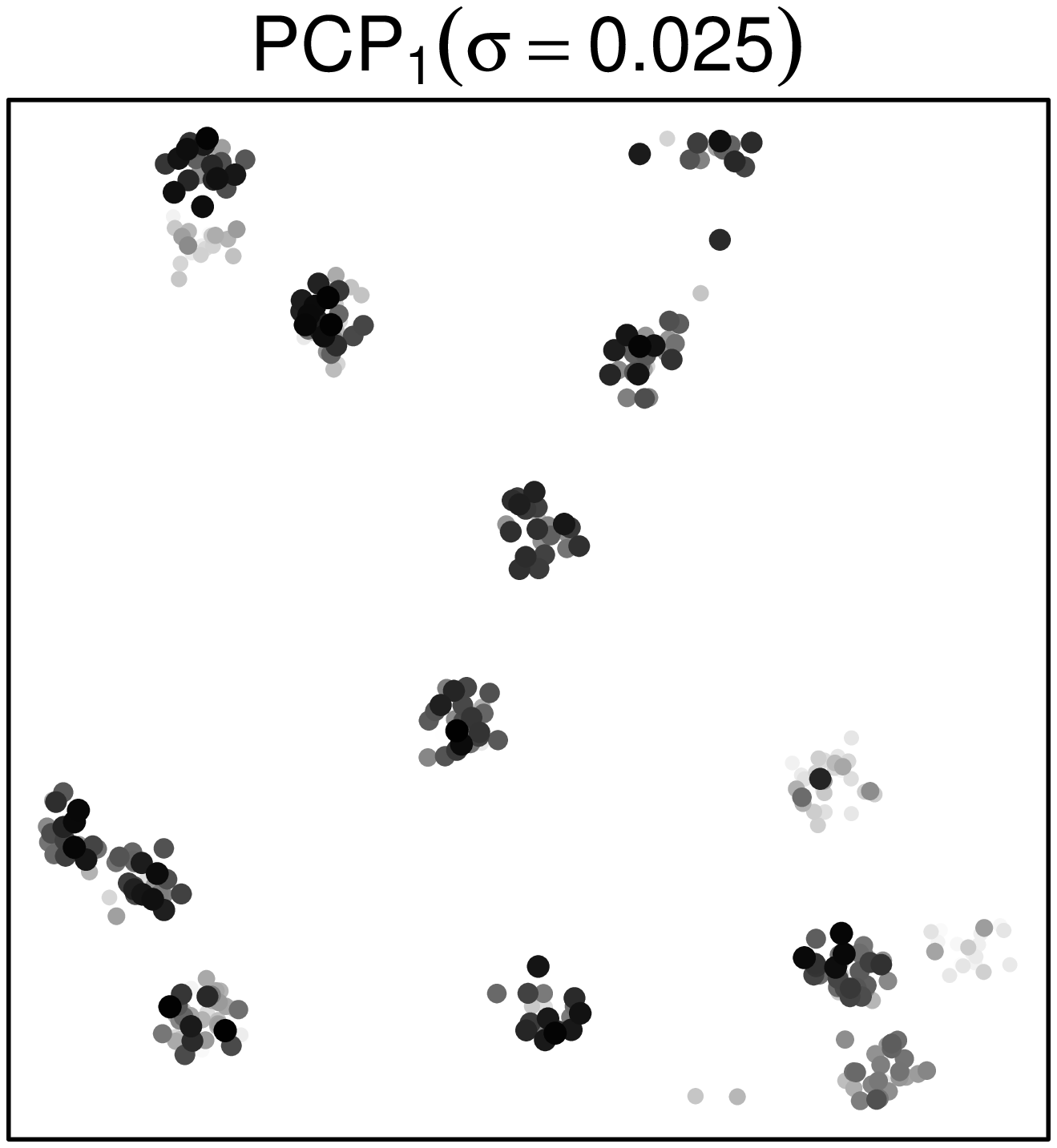}

  \includegraphics[width=3.75cm,height=3.75cm]{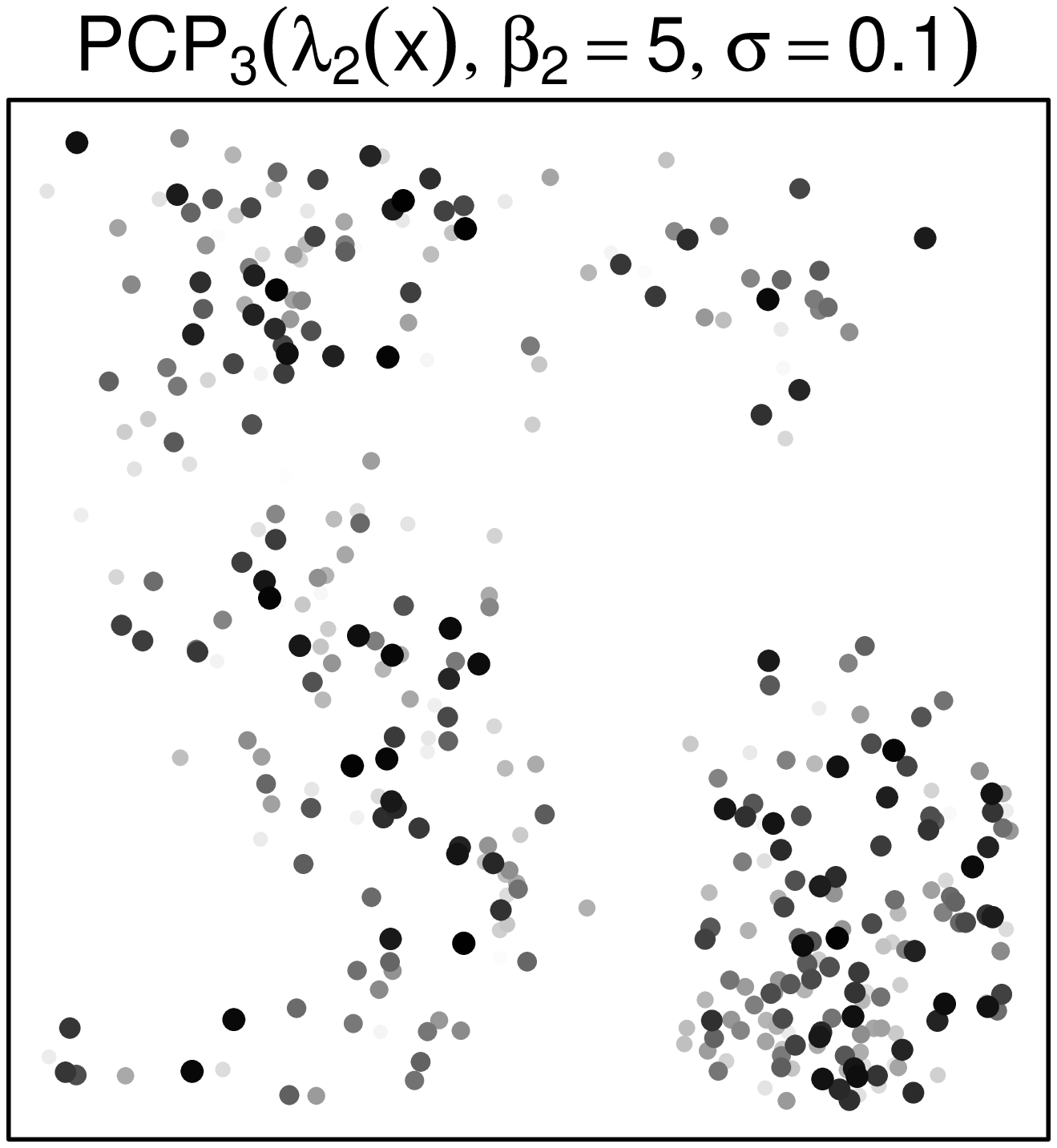}
    \includegraphics[width=3.75cm,height=3.75cm]{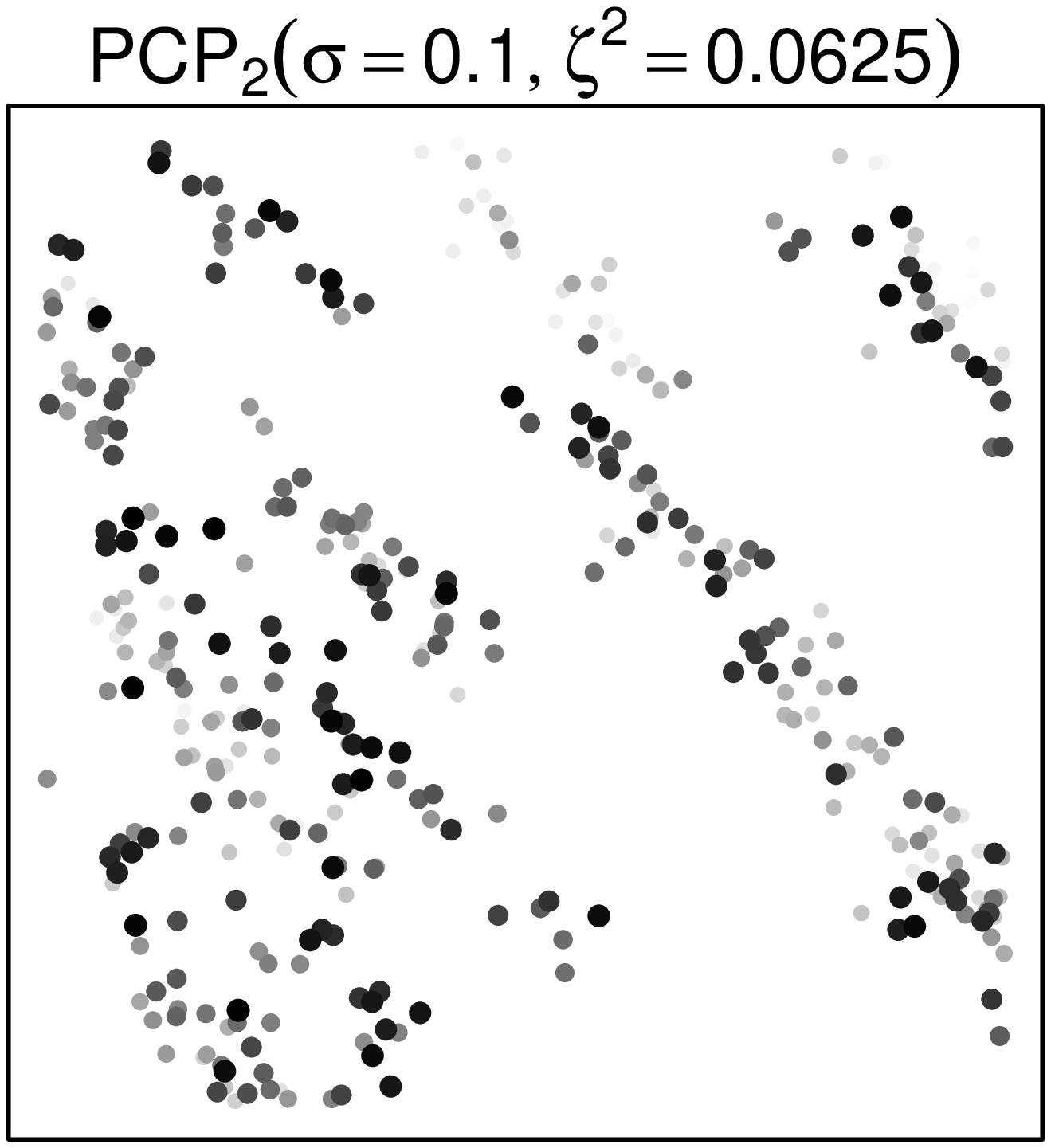}
      \includegraphics[width=3.75cm,height=3.75cm]{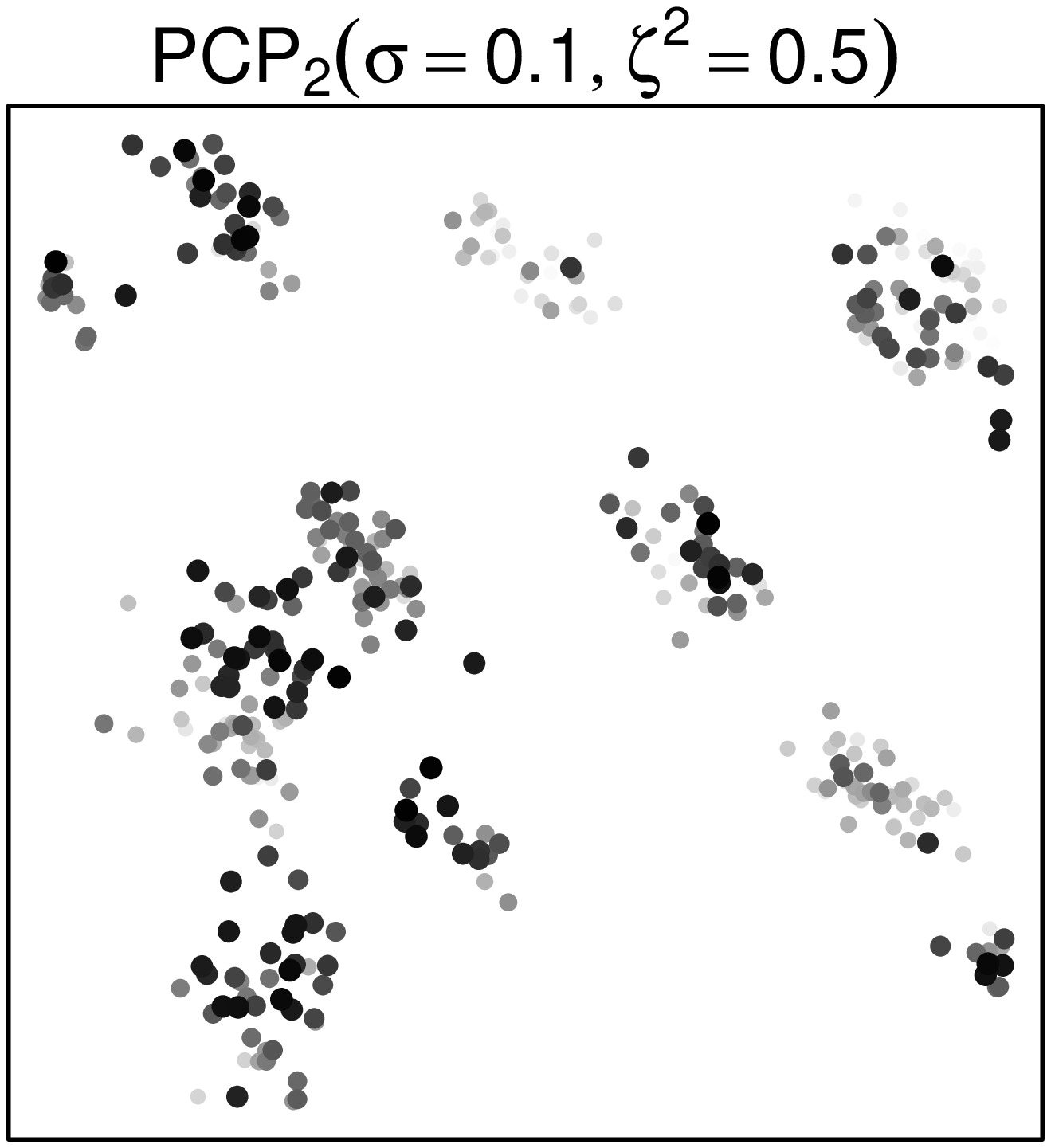}
 \end{center}
  \caption[]{\protect\parbox[t]{5in}{Realisations of the inhomogeneous Poisson et Poisson cluster processes. Dark and large dots correspond to recent events in time.}}
\label{fig:pp}
\end{figure}

For the simulation study, we set $S \times T = [0,1]^3$ and simulate processes with expected number of points $\bE [N(S \times T)] = n = 375$. We generate $N_{sim}=1000$ realisations of
\begin{itemize}
\item [(i)] $HPP(\lambda)$: an homogeneous Poisson process with intensity $\lambda=375$,

\item [(ii)] $IPP(\lambda(x);\beta)$: inhomogeneous Poisson processes with intensity function
\begin{equation}\label{eq:lambda}
    \lambda_1(x) = \lambda_1(s,t) = \frac{n \beta_1^3}{(e^\beta-1)^2(1-e^{-\beta_1})}\exp(\beta_1(s_x+s_y-t)),
\end{equation}
with $s=(s_x,s_y)$ and $\beta_1 \in \lce 1,2 \rce$,
or
\begin{equation}\label{eq:lambda2}
    \lambda_2(x) =  \frac{n \left( 1.25 + \cos(\beta_2 s_x +0.25) \right)
     \left( 1.25 + \cos(\beta_2 t +0.25) \right)}{\left(1.25 + \frac{1}{\beta_2} \left( \sin(\beta_2 + 0.25) -\sin(\beta_2) \right) \right)^2} ,
\end{equation}
with $\beta_2 \in \lce 3, 5 \rce$,

 \item [(iii)] $PCP_1(\nu,\sigma,\alpha,m_c)$: stationary Poisson cluster processes with intensity of parents $\nu=25$. The spatial distribution of the offspring is a zero-mean bivariate isotropic normal distribution
$\phi^{(2)}_{\sigma^2}$ with standard deviation $\sigma \in \lce 0.025, 0.05, 0.1, 0.15, 0.2 \rce$ and the temporal distribution is exponential with rate $\alpha=0.2$ and denoted ${\cal E}_\alpha(t)$. The expected number of offspring per parent follows a Poisson distribution with mean $m_c=15$.

This Poisson cluster process is an interpretation of a spatio-temporal shot-noise Cox process \citep{moller2010,moller2003b} with residual process $R(x)=\frac{1}{\nu} \sum_{y \in \Phi} \varphi(x-y)$,
where $\Phi$ is a stationary Poisson process in $\bR^2 \times \bR$ with intensity $\nu$ and $\varphi$ is the density function $\varphi(u,v)=\phi^{(2)}_{\sigma^2}(\|u\|){\cal E}_\alpha(v)$.
For such a process, we have $g(u,v)=1 + \frac{1}{\nu} \varphi \ast \tilde\varphi(u,v)$, where $\ast$ denotes convolution and $\tilde\varphi(u,v)=\varphi(-u,-v)$. Thus,
$$  g(u,v) =
    1 + \frac{\alpha}{8 \pi \sigma^2 \nu} \exp \left( - \frac{\|u\|^2}{4 \sigma^2} - \alpha v \right),$$
and
$$K(u,v) = 2 \pi u^2 v + \frac{1}{2 \nu}(1-\exp(-\alpha v)) \left(1-\exp \left(-\frac{u^2}{4 \sigma^2} \right) \right).$$
\item [(iv)] $PCP_2(\nu,\sigma,\alpha,m_c,\zeta^2,\theta,\omega)$: geometric anisotropic Poisson cluster processes, as defined in (iii) but with
$g(u,v)=g_0\left(\sqrt{u \Sigma^{-1} u^t},v \right)$,
where $u \in \bR^2$ is a row vector with transpose $u^t$, the function $g_0$ is such that $g$ is locally integrable \citep[see][]{moller2012b},
$\Sigma$ is a $2 \times 2$ symmetric positive definite matrix of the form $\Sigma=\omega^2 U_\theta \text{diag}(1,\zeta^2) U_\theta^t$ with $U_\theta = \left(
\begin{array}{cc}
  \cos(\theta) & -\sin(\theta) \\
  \sin(\theta) & \cos(\theta)
\end{array}
\right)$.
     Here $\theta=\pi/4$, $\omega=4$, and the anisotropy factor is $\zeta \in \lce 0.25,0.5,\sqrt{0.5},\sqrt{0.75},1 \rce$.

\item [(v)] $PCP_3(\lambda(x),\beta,\sigma,\alpha,m_c)$: non-stationary Poisson cluster processes, with intensity of parents defined by Equation~(\ref{eq:lambda}) or (\ref{eq:lambda2}) and clusters as defined in (iii).
\end{itemize}
Figure~\ref{fig:pp} shows static display of a realisation of the inhomogeneous Poisson process $IPP(\lambda_1(x),$ $\beta_1=2)$ (top left) and $IPP(\lambda_2(x),\beta_2=5)$ (middle top), the stationary Poisson cluster process with $\sigma=0.025$ (top right), the non-stationary Poisson cluster process $PCP_3(\lambda_2(x),\beta_2=5,\sigma=0.1)$ (bottom left) and
the anisotropic Poisson cluster process with $\sigma=0.1$  and $\zeta^2=0.0625$ (middle bottom) and $\zeta^2=0.5$ (bottom right).
The time is treated as a quantitative mark attached to each location, and the locations are plotted with the size and colour of the plotting symbol determined by the value of the mark. Dark and large dots correspond to recent events in time. Further forms of static and dynamic display of spatio-temporal point patterns can be found in \cite{gabriel2013}.

\subsection{Results}

The pair correlation function estimation is performed on the $N_{sim}=1000$ realised data sets of each point process. The STIK-function is only estimated for the processes (i)-(iii) and (v).
In both cases we use the different edge correction methods; we denote by $\widehat{K}_I$ and $\widehat{g}_I$, $\widehat{K}_B$ and $\widehat{g}_B$, $\widehat{K}_{MB}$ and $\widehat{g}_{MB}$, $\widehat{K}_T$ and $\widehat{g}_T$,
the estimators of the STIK-function and of the pair correlation function when using the isotropic, the border, the modified border, the translation edge correction factor respectively and $\widehat{K}_N$ and $\widehat{g}_N$ the estimators in the case of no edge correction.
In this section, we assume that the intensity is known to put ahead the influence of edge correction methods.

\subsection{Performance of the estimators relative to edge correction methods}

Performance is first measured through the empirical bias and mean squared error (MSE) which are computed for each spatial and temporal distances $u$ and $v$ in the sequence starting from 0.01 to 0.25 by increment of 0.01.

\begin{figure}[h]

 \begin{center}
 \vspace{-.5cm}
  \includegraphics[width=11.8cm,height=3cm]{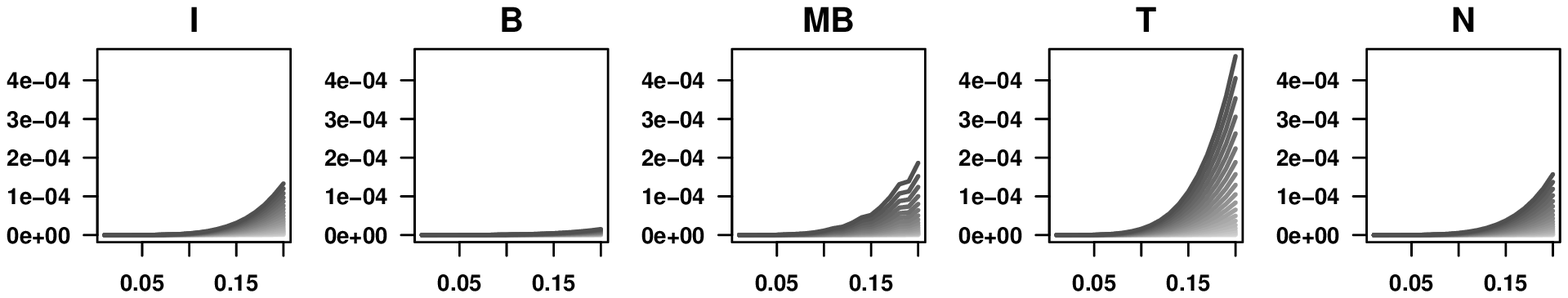}

\vspace{-.5cm}
 \includegraphics[width=11.8cm,height=3cm]{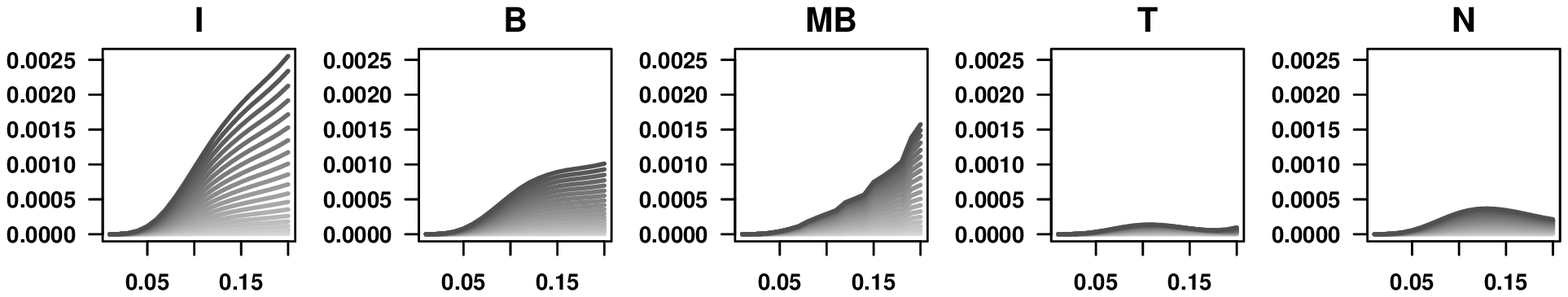}

\vspace{-.5cm}
  \includegraphics[width=11.8cm,height=3cm]{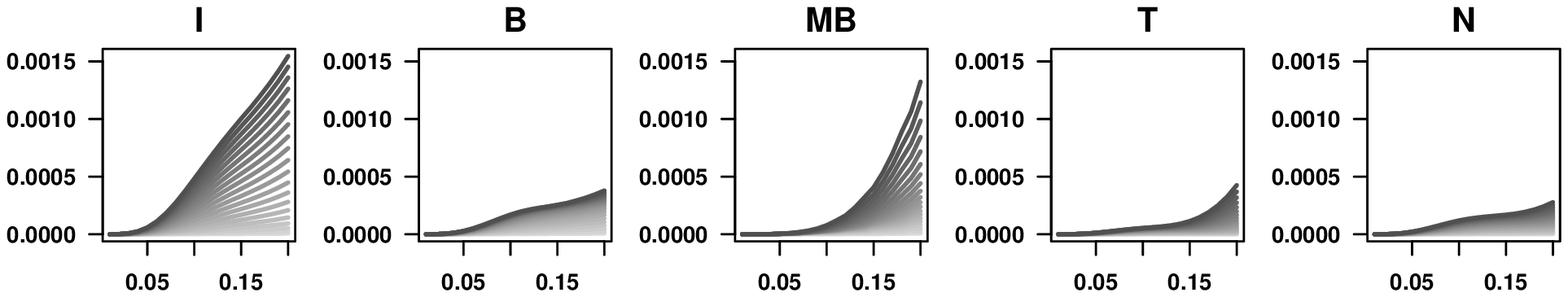}

  \vspace{-.5cm}
   \includegraphics[width=11.8cm,height=3cm]{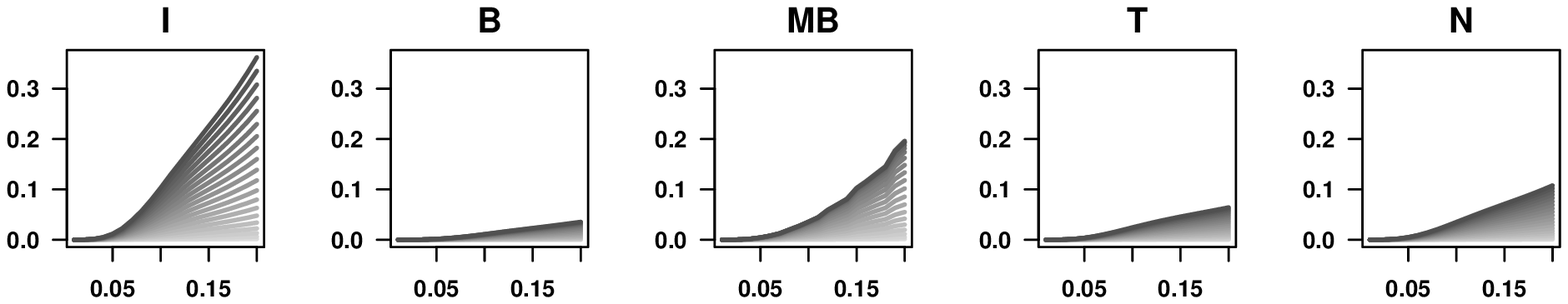}
 \end{center}

   \vspace{-.5cm}
     \caption[]{\protect\parbox[t]{5in}{Mean squared error of the STIK-function from $IPP(\lambda_1(x),\beta_1=1)$ (first row),
$PCP_1(\sigma=0.05)$ (second row), $PCP_3(\lambda_1(x),\beta_1=1,\sigma=0.05)$ (third row) and $PCP_3(\lambda_2(x),\beta_2=3,\sigma=0.05)$ (fourth row). Spatial distances $u$ are in abscissa; temporal distances are in grey: the darker, the greater, with the same range of values than $u$.}}
\label{fig:mseSTIK}
\end{figure}

\medskip
All estimators of the STIK-function have very small empirical bias and MSE at small spatial and temporal distances (less than 0.1) and their differences are observed for larger distances. Figure~\ref{fig:mseSTIK} shows the MSE of the STIK-function estimated when using the isotropic (I), border (B), modified-border (MB), translation (T), none (N) edge correction methods. Spatial distances $u$ are in abscissa, while temporal distances $v$ are in grey: the darker, the greater, with the same range of values than $u$.
For the Poisson processes (both homogeneous and inhomogeneous), $\widehat{K}_T$, $\widehat{K}_{MB}$ and $\widehat{K}_N$ have negative bias, while $\widehat{K}_I$ has a positive one. The bias of $\widehat{K}_B$ is almost constant. The MSE of $\widehat{K_T}$ has larger values than the others for large values of $u$ and $v$. This is illustrated for $IPP(\lambda_1(x),\beta_1=1)$ on the first row of Figure~\ref{fig:mseSTIK}.
Nevertheless, in all cases bias and MSE have negligible values.
Values of the bias and MSE slightly increase for the stationary Poisson cluster processes $PCP_1(\sigma)$ and particularly when the range of spatial clustering, $\sigma$, is small. For such processes, the bias of $\widehat{K}_B$ increases with $u$ and $v$,
$\widehat{K}_I$ and $\widehat{K}_T$ perform less well and better, respectively, than the others. The second row of Figure~\ref{fig:mseSTIK} shows the MSE of all estimators for the process $PCP_1(\sigma=0.05)$.
For the non-stationary Poisson cluster processes $PCP_3(\lambda_1(x),\beta_1,\sigma)$, bias and MSE of all estimators are similar to the stationary case. For $PCP_3(\lambda_2(x),\beta_2,\sigma)$, all estimators have positive bias. Bias and MSE increase with $\sigma$ and particularly for $\widehat{K}_{I}$ and $\widehat{K}_{MB}$.
The last two rows of Figure~\ref{fig:mseSTIK} illustrate the MSE of the estimators for the processes $PCP_3(\lambda_1(x),\beta_1=1,\sigma=0.05)$ (third row) and $PCP_3(\lambda_2(x),\beta_2=3,\sigma=0.05)$ (last row).

All estimators of the pair correlation function have small empirical bias and MSE with the exception
of the case of very small values of $u$ as the estimators contain the term ``$1/u$".  Bias at very small spatial distances could be reduced by the reflection technique \citep{doguwa1990}.
For the Poisson processes $HPP$ and $PCP_1(\lambda(x),\beta)$, we get similar comments than from the STIK-function, except that bias and MSE of $\widehat{g}_{MB}$ are almost constant with $u$ and that $\widehat{g}_{I}$ has a large variability at small distances, which increase with the inhomogeneity of $\lambda(x)$.
For the Poisson cluster processes $PCP_1(\sigma)$, $PCP_2(\sigma,\zeta^2)$ and $PCP_3(\lambda(x),\beta,\sigma)$, smaller the spatial dispersion of offspring $\sigma$ and stronger the anisotropy (i.e. greater $\zeta^2$), greater the bias and the MSE at small $u$. This is particularly so for $\widehat{g}_I$, $\widehat{g}_B$  and $\widehat{g}_{MB}$.
The MSE of all estimators of the pair correlation function are illustrated in Figure~\ref{fig:msePCF} (Appendix 2) for the same processes as in Figure~\ref{fig:mseSTIK} and for
$PCP_2(\sigma=0.2, \zeta^2=0.0625)$ and $PCP_2(\sigma=0.2, \zeta^2=0.5)$.

\medskip

The relative performance between the STIK-function and the pair correlation function does not only depend on the point process; it also depends on the edge correction method used to estimated the second-order characteristics.
We compare the performance of $\widehat{K}$ and $\widehat{g}$ from the absolute relative error (ARE), $|\widehat{T}(u,v) - T(u,v)|/T(u,v)$, evaluated for all $(u,v)$ and for each edge correction factor, with $T= K$ or $g$.
For the Poisson processes, the STIK-function perform better than the pair correlation function whatever the edge correction factor, except for the border and modified-border methods in the case of $IPP(\lambda_2(x),\beta_2=5)$ for which we get comparable performances of the two estimators.
For the Poisson cluster processes $PCP_1(\sigma)$ with $\sigma<0.15$, $\widehat{g}_I$ and $\widehat{g}_B$ have smaller ARE than $\widehat{K}_I$ and $\widehat{K}_B$ for all distances, but not for the others at large distances. When $\sigma>0.15$, this phenomenon is stronger and $\widehat{K}_T$ can globally be better than $\widehat{g}_T$. We get similar results for $PCP_3(\lambda_1(x),\beta_1,\sigma)$. For $PCP_3(\lambda_2(x),\beta_2,\sigma)$, $\widehat{g}$ have to be preferred whatever the edge correction method used in the estimation.

\medskip

To avoid problems associated with multiple testing at different spatial and temporal scales, we also compute integral deviation measures
\begin{equation}\label{eq:dev}
D=\int_0^{v'} \int_0^{u'} ( \widehat T(u,v) - T(u,v) )^2 \dd u \dd v,
\end{equation}
where the summary statistic $T$ is either $K$ or $g$.
We then compute the relative efficiency of the overall $D$ defined by
$100 \times \min_k V_k / V_k$, where $V_k$ is the empirical variance
 $V_k= \dfrac{1}{N_{sim}} \sum_{i=1}^{N_{sim}} \left(D_{ik} - \overline{D_{k}} \right)^2$, $k \in \lce I, B, MB, T, N \rce $, with $D_{ik}$ the deviation measure evaluated from the estimation of $T$ for the $i$th simulation and the $k$th edge correction method
 and  $\overline{D_k}=\frac{1}{N_{sim}} \sum_{i=1}^{N_{sim}} D_{ik}$.
 Tables~\ref{tab:devedgeK} and \ref{tab:devedge} (see Appendix 2) give the relative efficiency of the overall $D$ defined from the STIK-function and the pair correlation function respectively. They are summarized in Table~\ref{tab:dev} which provides, for each point process, the best edge correction method, i.e. the one for which the relative efficiency is 100.  It shows that the border method is the most efficient for Poisson processes (all cases for $\widehat{K}$ and homogeneous or with weak inhomogeneity for $\widehat{g}$) and for the non-stationary Poisson cluster processes $PCP_3(\lambda_2(x),\beta_2)$. Then for inhomogeneous and/or clustered isotropic and anisotropic processes, the translation method is the most efficient.
 \begin{table}[h]
 \centering
\caption[]{\protect\parbox[t]{5in}{Most efficient edge correction method according to the point process.}}
  \begin{tabular}{|l|cc|}
      \hline
& $\widehat{K}$ & $\widehat{g}$ \\
              \hline
$HPP(\lambda)$ & Border & Border \\
$IPP(\lambda(x),\beta)$ & Border & Translation \\
$PCP_1(\sigma)$ & Translation & Translation \\
$PCP_2(\sigma, \zeta^2)$ & - & Translation \\
$PCP_3(\lambda_1(x),\beta_1=1, \sigma)$ & Translation & Translation \\
$PCP_3(\lambda_1(x),\beta_1=2, \sigma)$ & Modified Border & Translation \\
$PCP_3(\lambda_2(x),\beta_2, \sigma)$ & Border & Border \\
      \hline
    \end{tabular}
\label{tab:dev}
\end{table}

\subsection{Clustering detection relative to edge correction methods}
\label{subsec:clust}

Let us now look at the ability of detecting clustering relative to the edge correction method. The pair correlation function is evaluated from the realisations of the Poisson cluster processes $PCP_1(\sigma)$ and $PCP_3(\lambda(x),\beta,\sigma)$.
We say that there is a clustering tendency when $\widehat g$ is greater than the upper envelope of the pair correlation function evaluated under the null hypothesis (absence of clustering), $g_{up}(u,v)=\max_{i} \widehat{g_{ik}}(u,v)$, i.e. evaluated for the Poisson processes $HPP(\lambda)$ and $IPP(\lambda(x),\beta)$. Thus, for all distances $u$ and $v$ and for each edge correction method (subscript $k$), we compute the probability of detecting clustering:
$p_k(u,v) = \frac{1}{N_{sim}} \sum_{i=1}^{N_{sim}} {\bf 1}_{\lce \widehat{g_{ik}}(u,v) > g_{up}(u,v) \rce}$.
The envelopes are built from $1000$ simulations, thus leading to reasonable values for the type I error probabilities (between 0.02 and 0.04 for $\widehat K$ and between 0.04 and 0.1 for $\widehat g$ according to the edge correction method); see \cite{loosmore2006} and \cite{grabarnik2011} for a guidance on the use of envelop tests and deviation tests.
The probability of detecting clustering is expected to be maximum for spatial distances $u$ less than $\approx 2 \sigma$ corresponding to the size of the clusters defined by our Poisson cluster processes. Because of the range of spatial distances used to estimate the pair correlation (between 0.01 and 0.25), we restrict our analysis to $\sigma \in \lce 0.025, 0.05, 0.1 \rce$. Results from the STIK-function are omitted for conciseness.

\begin{figure}[h]
  \begin{center}
\parbox{4.5cm}{\sloppy
\begin{center}
{\tiny $PCP_1(\sigma=0.025)$}

 \vspace{-4mm}
  \includegraphics[width=3.8cm,height=2.5cm]{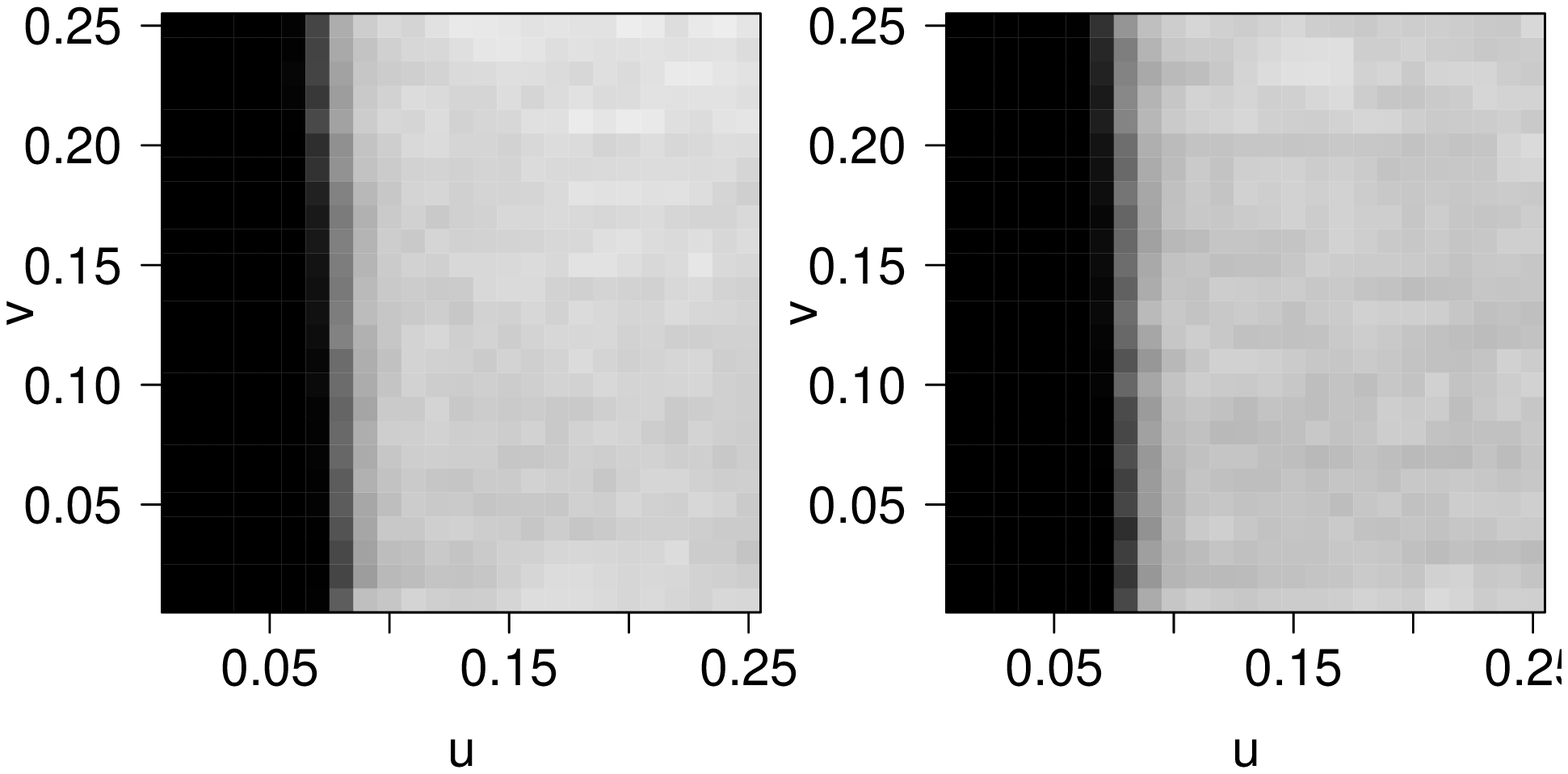}

   \vspace{-3mm}
{\tiny $PCP_3(\lambda_1(x),\beta_1=1,\sigma=0.025)$}

 \vspace{-4mm}
  \includegraphics[width=3.8cm,height=2.5cm]{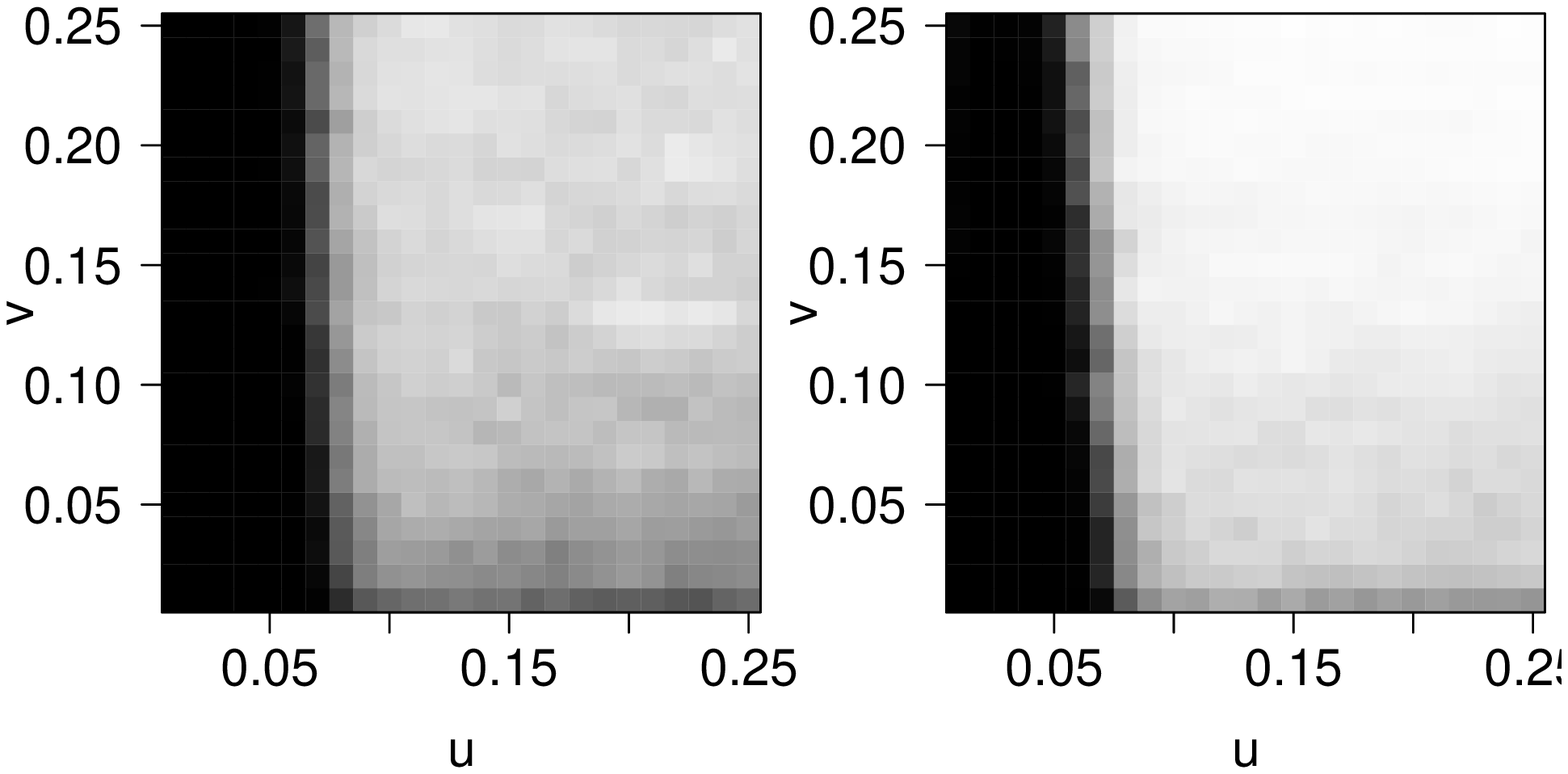}

   \vspace{-3mm}
  {\tiny $PCP_3(\lambda_2(x),\beta_2=3,\sigma=0.025)$}

 \vspace{-4mm}
  \includegraphics[width=3.8cm,height=2.5cm]{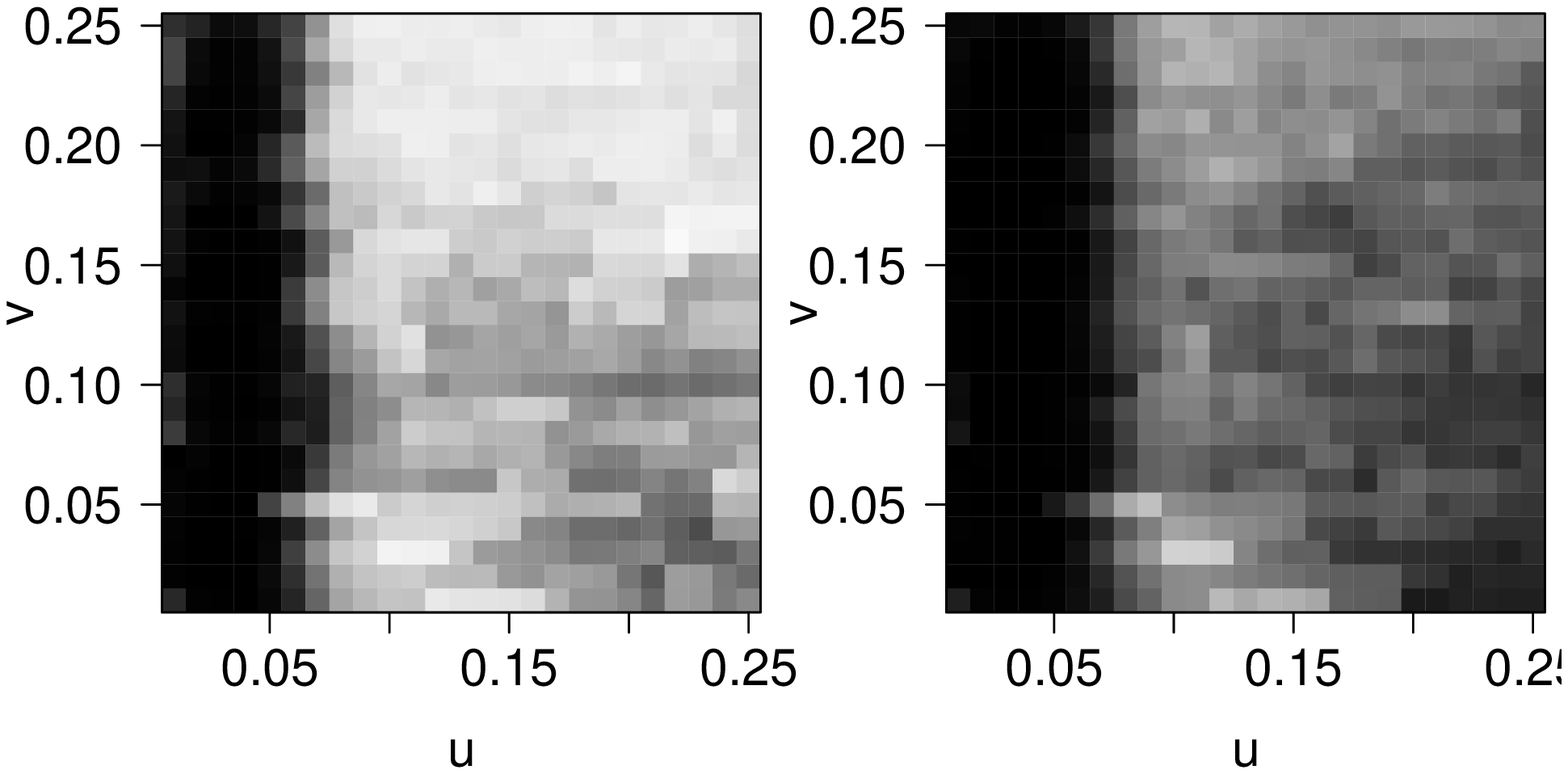}
\end{center}
}
\parbox{4.5cm}{\sloppy
\begin{center}
{\tiny $PCP_1(\sigma=0.05)$}

 \vspace{-4mm}
  \includegraphics[width=3.8cm,height=2.5cm]{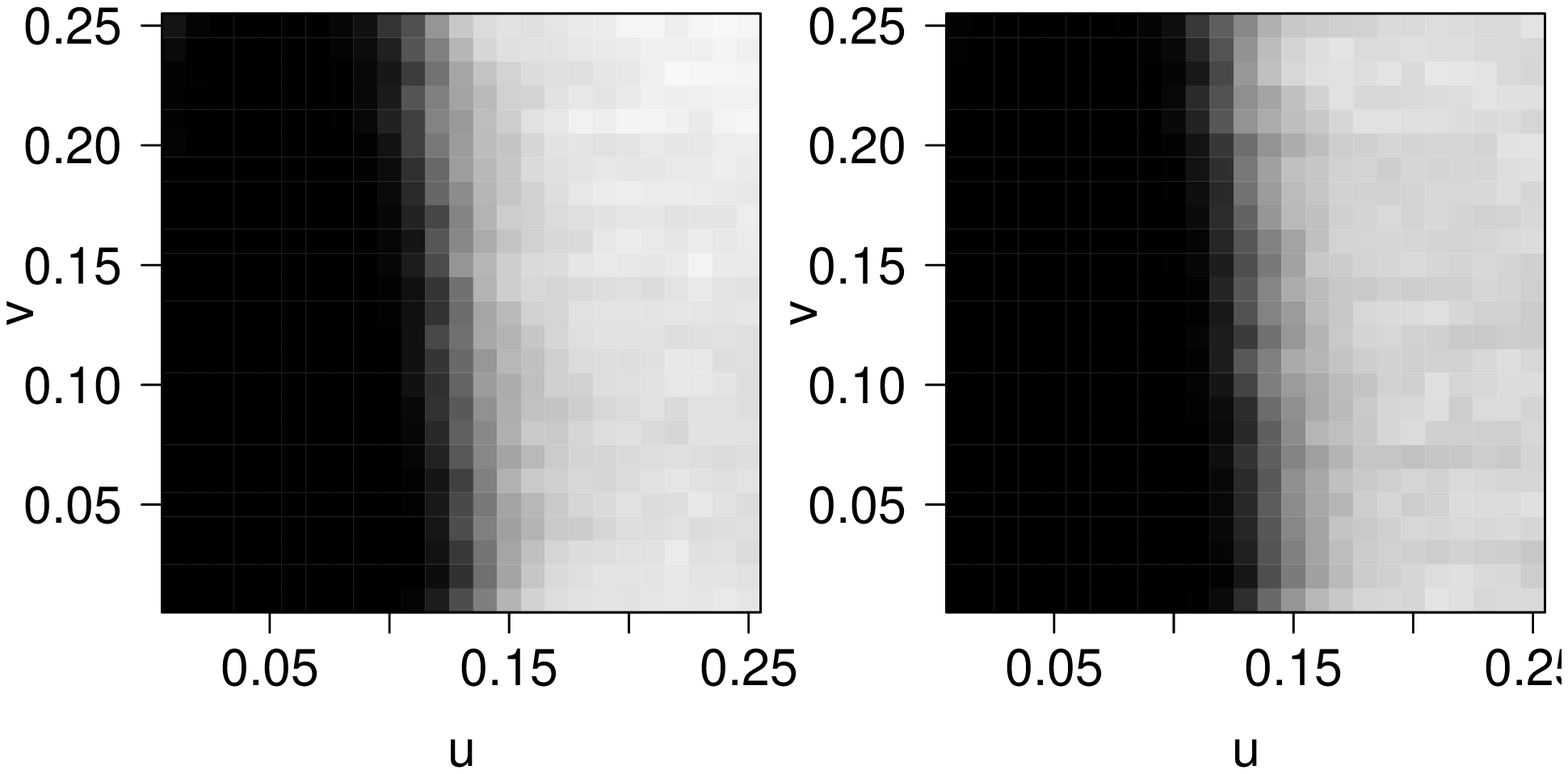}

 \vspace{-3mm}
{\tiny $PCP_3(\lambda_1(x),\beta_1=2,\sigma=0.025)$}

 \vspace{-4mm}
  \includegraphics[width=3.8cm,height=2.5cm]{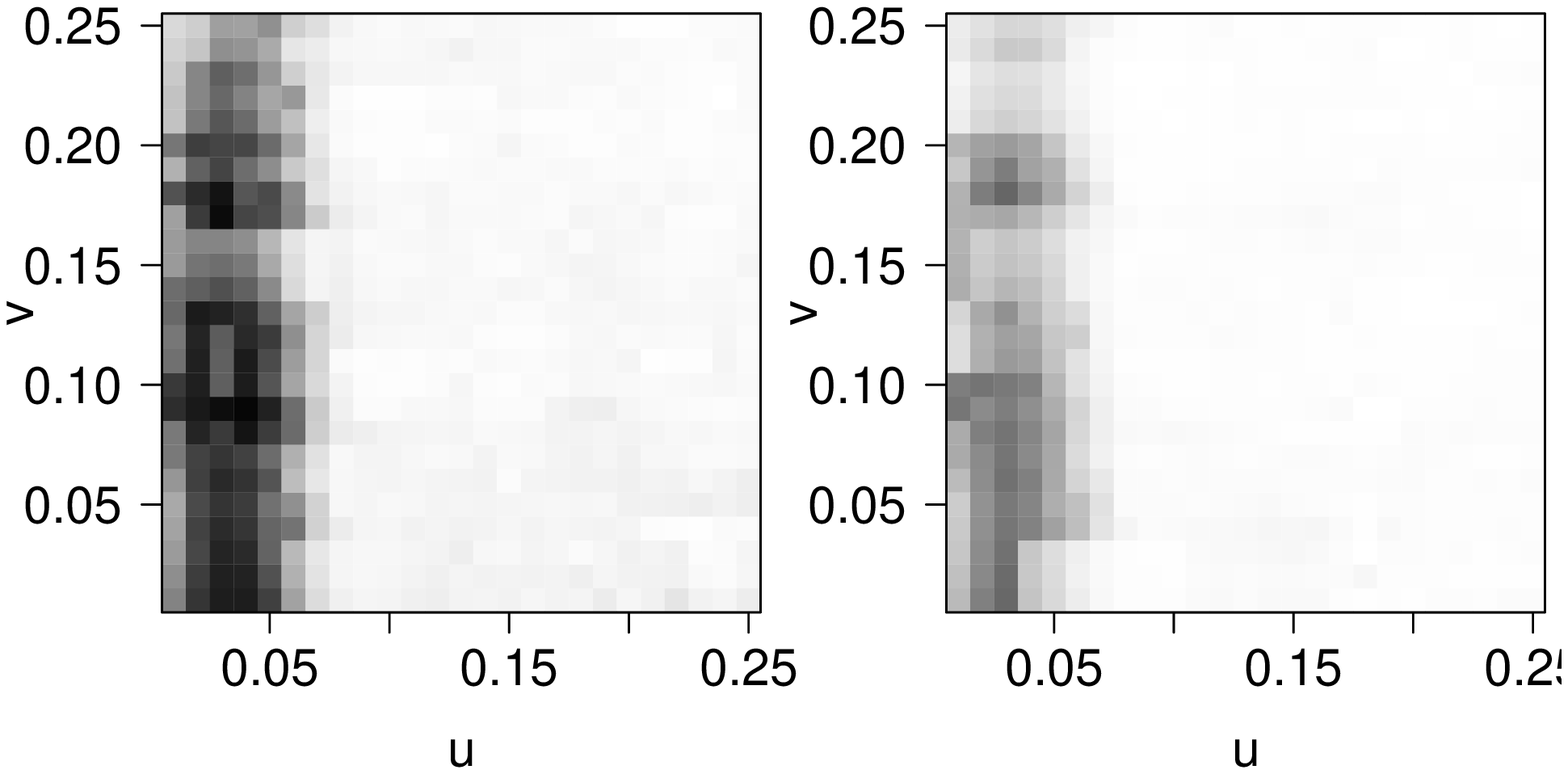}

 \vspace{-3mm}
  {\tiny $PCP_3(\lambda_2(x),\beta_2=5,\sigma=0.025)$}

 \vspace{-4mm}
  \includegraphics[width=3.8cm,height=2.5cm]{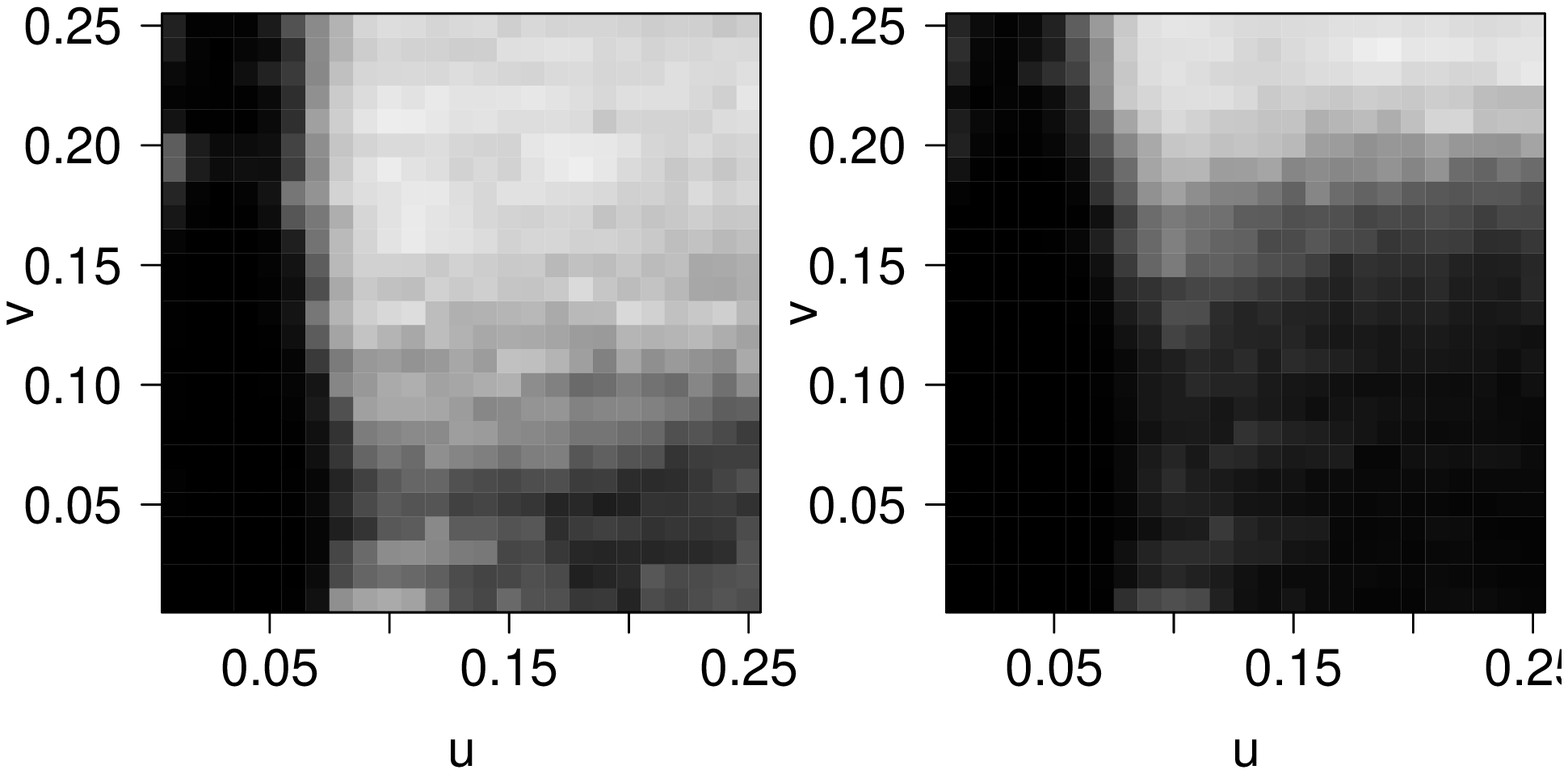}
\end{center}
}

\end{center}

\vspace{-.5cm}
  \caption[]{\protect\parbox[t]{5in}{Probabilities of detecting clustering using the border method (first and third columns) and translation method (second and fourth columns)
  for the Poisson cluster processes $PCP_1(\sigma)$ and $PCP_3(\lambda(x),\beta,\sigma)$.}}
  \label{fig:prob}
\end{figure}

For the stationary Poisson cluster processes $PCP_1(\sigma)$, the probability is maximum for spatial distances up to $0.07$, $0.12$ and $0.18$ when $\sigma = 0.025$, 0.05 and 0.1 respectively. This is observed for the border, modified-border and translation methods and when no edge correction is done. The isotropic edge correction factor leads to less powerful results as $\sigma$ increases.
We get similar results for $PCP_3(\lambda_1(x),\beta_1=1,\sigma)$,  than for $PCP_1(\sigma=0.025)$. However, for $\widehat{g}_{MB}$, $\widehat{g}_T$ and $\widehat{g}_N$ the probability of detecting clustering decreases both when $\sigma$ and $v$ increase.
This is also the case for $PCP_3(\lambda_1(x),\beta_1=2,\sigma)$, but the probabilities are weaker for all $u$ and $v$ because, whatever the edge correction method used, $g$ is under-smoothed for this process. This can be explained by some overlapped clusters.
For the processes $PCP_3(\lambda_2,\beta_2,\sigma)$, $\widehat{g}_I$ is only able to detect clustering when $\sigma=0.025$. The estimators $\widehat{g}_{B}$, $\widehat{g}_{MB}$, $\widehat{g}_T$ and $\widehat{g}_N$ behave similarly than for $PCP_3(\lambda_1(x),\beta_1=1,\sigma)$, except that ``false detections'' appear for large values of $u$. This is straightly related to over-smoothed second-order characteristics for such processes.
To illustrate all these comments, the probability of detecting clustering of $\widehat{g}_B$ and $\widehat{g}_T$ are plotted in Figure~\ref{fig:prob} for the Poisson cluster processes $PCP_1(\sigma=0.025)$, $PCP_1(\sigma=0.05)$ and $PCP_3(\lambda(x),\beta,\sigma)$ with a spatial
  dispersion $\sigma=$ 0.025. Grey shading indicates values of the probability of detecting clustering. Black corresponds to one and white corresponds to 0.
\clearpage
Similarly, we obtain probabilities of detecting clustering for the overall statistic $D$. However, to avoid misleading results due to bias in the estimations, we replace $D$ defined in Equation~(\ref{eq:dev}) by
\begin{equation}\label{eq:dev2}
    \widetilde D_{ik}=\int_0^{v'} \int_0^{u'} ( \widehat g_{ik}(u,v) - \overline{\widehat g_k(u,v)} )^2 \dd u \dd v,
\end{equation}
for $k \in \lce I, B, MB, T, N\rce$, $i=1,\dots, N_{sim}$ and $\overline{\widehat g_k(u,v)}=\frac{1}{N_{sim}} \sum_{i=1}^{N_{sim}} \widehat g_{ik}(u,v)$.
Table~\ref{tab:power} shows these probabilities for various scales of spatial clustering and inhomogeneity. The performance of the estimators are good for the stationary Poisson cluster processes $PCP_1(\sigma)$. As soon as we introduce spatial inhomogeneity, $\widehat{g}_I$ is no longer able to detect clustering. The border method appear to provide the best results. For $PCP_3(\lambda(x),\beta,\sigma)$, the power of each estimator is related to the range of spatial clustering $\sigma$, it decreases as $\sigma$ increases, and to the model of inhomogeneity $\lambda(x)$.
\begin{table}[h]
\centering
\caption[]{\protect\parbox[t]{5in}{Power of the integral deviation measures $\widetilde D$ according to the isotropic (I), border (B), modified border (MB), translation (T) edge correction factor and without edge correction (N).}}
    \begin{tabular}{|l|ccccc|}
      \hline
        & \multicolumn{5}{|c|}{Edge correction factor} \\
             Point process & I & B & MB & T & N \\
              \hline
$PCP_1(\sigma=0.025)$ & 1.000 & 1.000 &     1.000 &    1.000& 1.000 \\
$PCP_1(\sigma=0.05)$  & 0.856 & 1.000 &     1.000 &    1.000& 1.000 \\
$PCP_1(\sigma=0.1)$   & 0.007 & 1.000 &     0.999 &    1.000& 1.000 \\
$PCP_3(\lambda_1(x),\beta_1=1,\sigma=0.025)$ & 0.667 & 1.000 &     0.961 &    1.000 & 1.000 \\
$PCP_3(\lambda_1(x),\beta_1=1,\sigma=0.05)$  & 0.208 & 1.000 &     0.625 &    0.947 & 0.927 \\
$PCP_3(\lambda_1(x),\beta_1=1,\sigma=0.1)$   & 0.091 & 0.968 &     0.338 &    0.570 &0.560 \\
$PCP_3(\lambda_1(x),\beta_1=2,\sigma=0.025)$ & 0.187 & 0.789 &    0.458 &    0.495 &0.411 \\
$PCP_3(\lambda_1(x),\beta_1=2,\sigma=0.05)$  & 0.077 & 0.216 &     0.129 &    0.129 & 0.112 \\
$PCP_3(\lambda_1(x),\beta_1=2,\sigma=0.1)$   & 0.042 & 0.089 &     0.079  &   0.073 & 0.069 \\
$PCP_3(\lambda_2(x),\beta_2=3,\sigma=0.025)$ & 0.226 & 1.000  &    1.000 &    1.000 &1.000 \\
$PCP_3(\lambda_2(x),\beta_2=3,\sigma=0.05)$  & 0.145 & 0.976  &    1.000 &    1.000 &1.000 \\
$PCP_3(\lambda_2(x),\beta_2=3,\sigma=0.1)$   & 0.089 & 0.776  &    1.000  &   1.000 &1.000 \\
$PCP_3(\lambda_2(x),\beta_2=5,\sigma=0.025)$ & 0.826 & 1.000  &    1.000  &   1.000 &1.000 \\
$PCP_3(\lambda_2(x),\beta_2=5,\sigma=0.05)$  & 0.476 & 1.000  &    1.000  &   1.000 &1.000 \\
$PCP_3(\lambda_2(x),\beta_2=5,\sigma=0.1)$   & 0.385 & 1.000  &    1.000  &   1.000 &1.000 \\
      \hline
    \end{tabular}
\label{tab:power}
\end{table}

\section{Influence of the intensity estimate}
\label{sec:3}

The estimators of the second-order characteristics defined in Section~\ref{subsec:STIK} are approximatively unbiased. However, they depend on the first-order intensity function $\lambda(x)$. As $\lambda(x)$ is not known in practice, it must be replaced by an estimate, $\hat \lambda(x)$, in Equations~(\ref{eq:STIKhat1}) and (\ref{eq:ghat}). There are several ways in estimating the intensity function; see for example \cite{cressie1993}, \cite{illian2008} or \cite{stoyan1994} for a review. The two mains are using a parametric model, often used when the point pattern suggest a theoretical form of $\lambda(x)$, or using kernel methods.
In this section, we illustrate by simulations the influence of $\hat \lambda(x)$ on the estimation of the pair correlation function.

\subsection{Estimation}

In the following, we assume the spatio-temporal separability of $\lambda(x)$. Thus, whatever the method used, the estimation of the intensity is treated completely separately for space and time.

\subsubsection{Parametric estimation}

We consider an exponential model for the temporal intensity function and a log linear function of spatial coordinates for the spatial intensity \citep{liu2007,stoyan1994}: $\lambda(s;\beta) = \exp(\beta^T Z(s))$, where $Z(s)$ is a vector of polynomials of coordinates $s_x,s_y$ of point location $s$. Polynomials are of order one and two:
\begin{eqnarray*}
\text{Model 1: } \lambda(s;\beta) & = & \exp(\beta_0 + \beta_1 s_x + \beta_2 s_y). \\
\text{Model 2: } \lambda(s;\beta) & = & \exp(\beta_0 + \beta_1 s_x + \beta_2 s_y + \beta_3 s_x^2 + \beta_4 s_x s_y + \beta_5 s^2_y).
\end{eqnarray*}
Parameters are estimated using likelihood methods. The log-likelihood for $\beta$ is
$$\log L(\beta,{\bf s}) = \sum_{i=1}^{n} \log (\lambda(s_i;\beta)) - \int_S \lambda(u;\beta) \dd u.$$ We use the Berman-Turner algorithm \citep{berman1992}
for finding the maximum likelihood estimate of $\beta$.

\subsubsection{Kernel estimation}
\label{subsec:kern}

The temporal intensity function is estimated using a gaussian kernel estimator, with bandwidth $0.9 n^{-1/5} \min \lce \sigma,(Q_3-Q_1)/1.34 \rce$, see \cite{silverman1986}.
To estimate the spatial intensity function we use a quartic kernel estimator. The estimator takes the form
$$\hat \lambda(s) = \sum_{i=1}^n\frac{k_h(s-s_i)}{c_{S,h}(s_i)},$$
where $k_h(s)= h^{-2} k(s/h)$ and $c_{S,h}(s_i)=\int_S k_h(s-s_i) \dd s$ is the edge correction factor defined in \cite{diggle1985} ensuring that $\int_S  \hat \lambda(s) \dd s = n$. The bandwidth, $\hat h$,  minimises the mean squared error \citep{berman1989}.

\subsection{Simulations}

Simulations are based on the 1000 realisations of the inhomogeneous Poisson processes $IPP(\lambda(x),\beta)$  and of the Poisson cluster processes $PCP_3(\lambda(x),\beta,\sigma)$ described in Section~\ref{subsec:simu}. According to the results of Section 3, we use the translation method to correct edge effects and consider the pair correlation function.
For each simulation, we compute the pair correlation function according to the different estimates of the intensity function. In the following,
the estimator is indexed by $\lambda$ when using the true intensity, by $\hat \lambda_p$ when using a parametric estimate and by $\hat \lambda_k$ when using a kernel estimate of the intensity function.

\subsection{Performance of $\widehat g$ relative to $\hat \lambda$}

\begin{figure}[h]
  \begin{center}
  \parbox{9cm}{\sloppy
  \begin{center}
  {\tiny $IPP(\lambda_1(x),\beta_1=1)$}

    \vspace{-0.4cm}
\includegraphics[width=8cm,height=2.5cm]{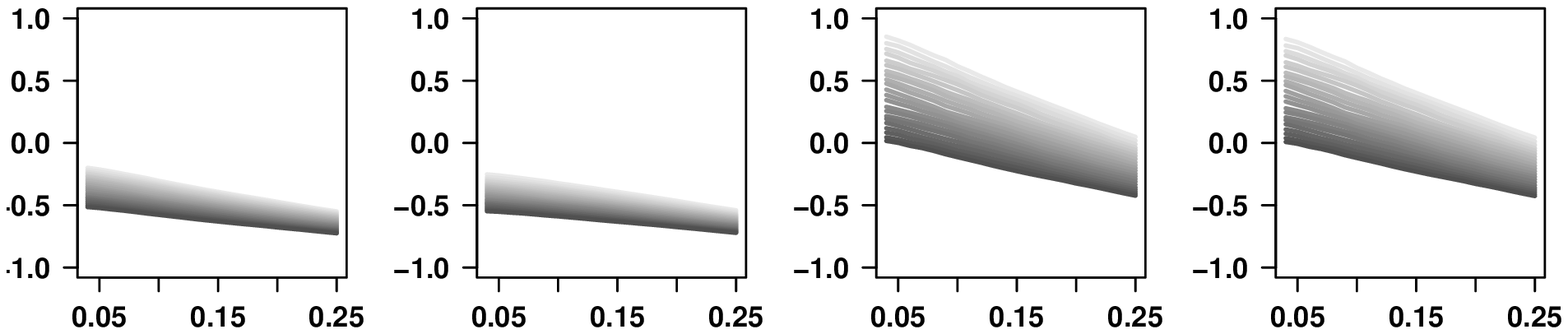}
  \end{center}
  }

  \vspace{-0.7cm}
  \parbox{9cm}{\sloppy
  \begin{center}
  {\tiny $IPP(\lambda_2(x),\beta_2=3)$}

    \vspace{-0.4cm}
\includegraphics[width=8cm,height=2.5cm]{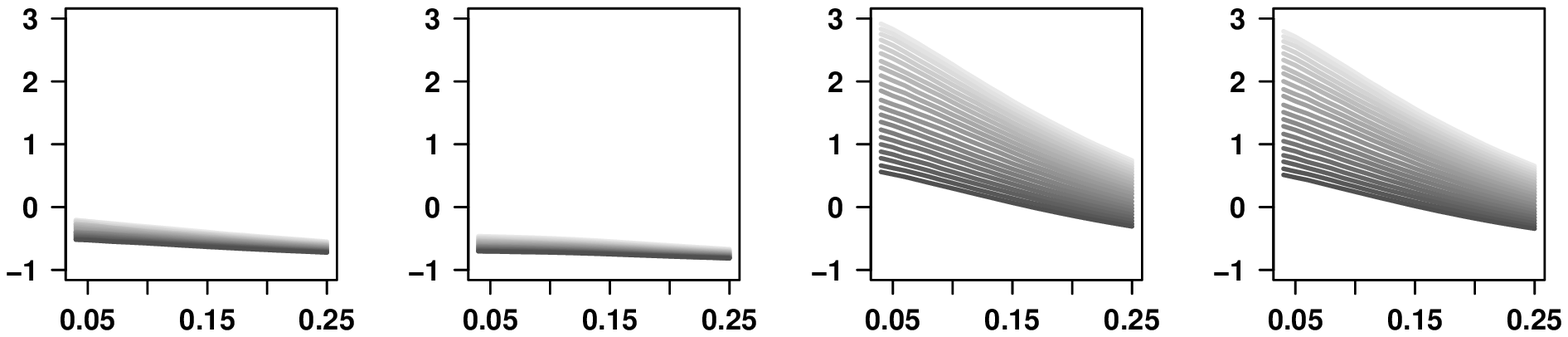}
  \end{center}
  }

    \vspace{-0.7cm}
    \parbox{9cm}{\sloppy
  \begin{center}
  {\tiny $PCP_3(\lambda_1(x),\beta_1=1,\sigma=0.05)$}

    \vspace{-0.4cm}
\includegraphics[width=8cm,height=2.5cm]{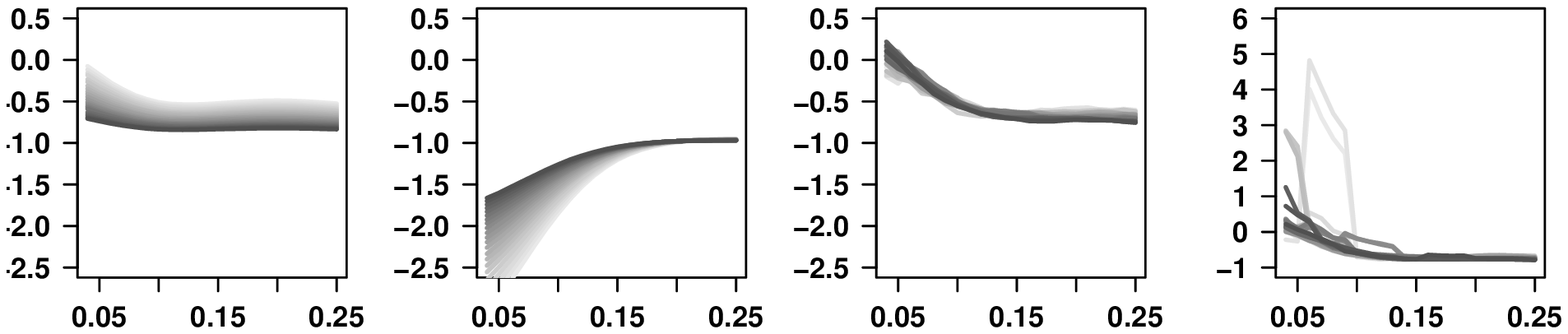}
  \end{center}
  }

    \vspace{-0.7cm}
    \parbox{9cm}{\sloppy
  \begin{center}
  {\tiny $PCP_3(\lambda_2(x),\beta_2=3,\sigma=0.05)$}

    \vspace{-0.4cm}
\includegraphics[width=8cm,height=2.5cm]{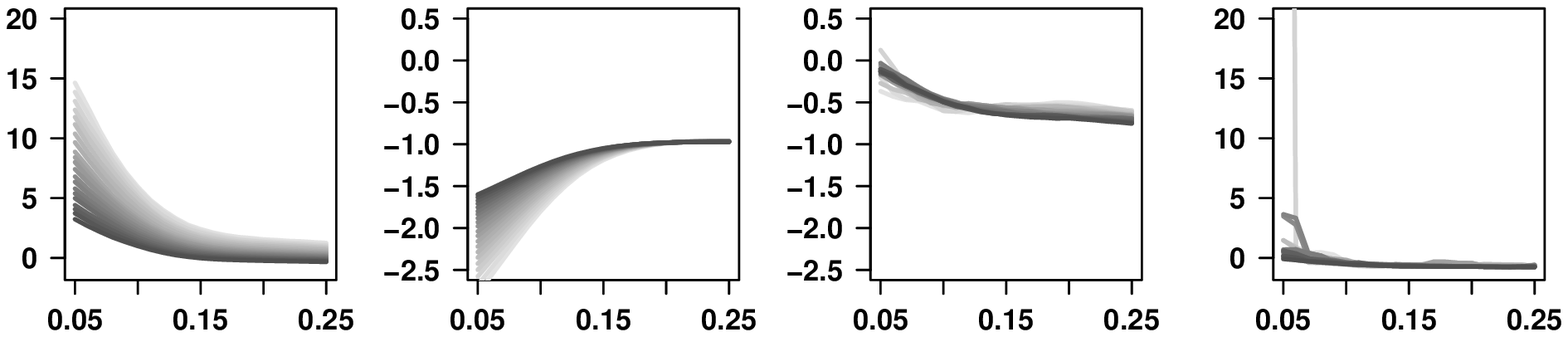}
  \end{center}
  }
  \end{center}

    \vspace{-0.7cm}
  \caption[]{\protect\parbox[t]{5in}{Bias of the pair correlation function estimated by using ${\lambda}$ (first column), $\widehat \lambda_k$ (second column), $\widehat \lambda_{p_1}$ (third column) and $\widehat \lambda_{p_2}$ (fourth column). Spatial distances $u$ are in abscissa, while temporal distances $v$ are in grey with the same values than $u$ (the darker, the greater).}}
  \label{fig:biaisPCF}
\end{figure}

We compute the pointwise bias and mean squared error of the pair correlation function estimator according to the different estimates of the intensity function.
For the inhomogeneous Poisson process $IPP(\lambda(x),\beta)$, $\widehat{g_\lambda}$ and $\widehat{g_{\hat \lambda_k}}$ have similar MSE and negative bias, both increasing with $u$ and $v$. Bias and MSE of $\widehat{g_{\hat \lambda_p}}$ are greater than for the two others, maybe due to over-smoothing, and much more greater for small values of $u$ and $v$ and as the inhomogeneity is stronger. The bias is positive for small $u$ and $v$ and becomes negative for large $u$ and $v$. This is illustrated in the first two rows of Figure~\ref{fig:biaisPCF}.
The last two rows of Figure~\ref{fig:biaisPCF} show results for the Poisson cluster processes. Model 2 implies very high bias at different ranges of spatial and temporal distances, thus leading to very high values of the corresponding MSE. For $PCP_3(\lambda_1(x),\beta_1,\sigma)$ the bias and MSE of $\widehat{g_\lambda}$ and  $\widehat{g_{\hat \lambda_p}}$ (Model 1) have similar orders of magnitude, except for small $u$, where they are greater when using parametric estimates of the intensity. This is not the case for $PCP_3(\lambda_2(x),\beta_2,\sigma)$.
The estimator $\widehat{g_{\hat \lambda_k}}$ is biased for the Poisson cluster processes.
This is mainly due to the fact that we estimated non parametrically both the intensity function and the pair correlation function from the same observed point pattern. Without any other assumption or information about the nature of the underlying point process, we cannot in general make an empirical distinction between first-order and second-order effects and thus cannot distinguish the contributions due to spatial inhomogeneity or spatial interacting phenomena.
Here, the bandwidths $\hat h$ selected for spatial kernel smoothing (see Section~\ref{subsec:kern}) are closed to the spatial dispersion $\sigma$. This leads to over-smoothed intensity functions and thus to severely biased values of the pair correlation function.

\subsection{Clustering detection relative to $\hat \lambda$}

\begin{figure}[h]
  \begin{center}
  {\tiny $PCP_3(\lambda_1(x),\beta_1=1,\sigma)$}

    \vspace{-1cm}
\includegraphics[width=7.5cm,height=4.2cm]{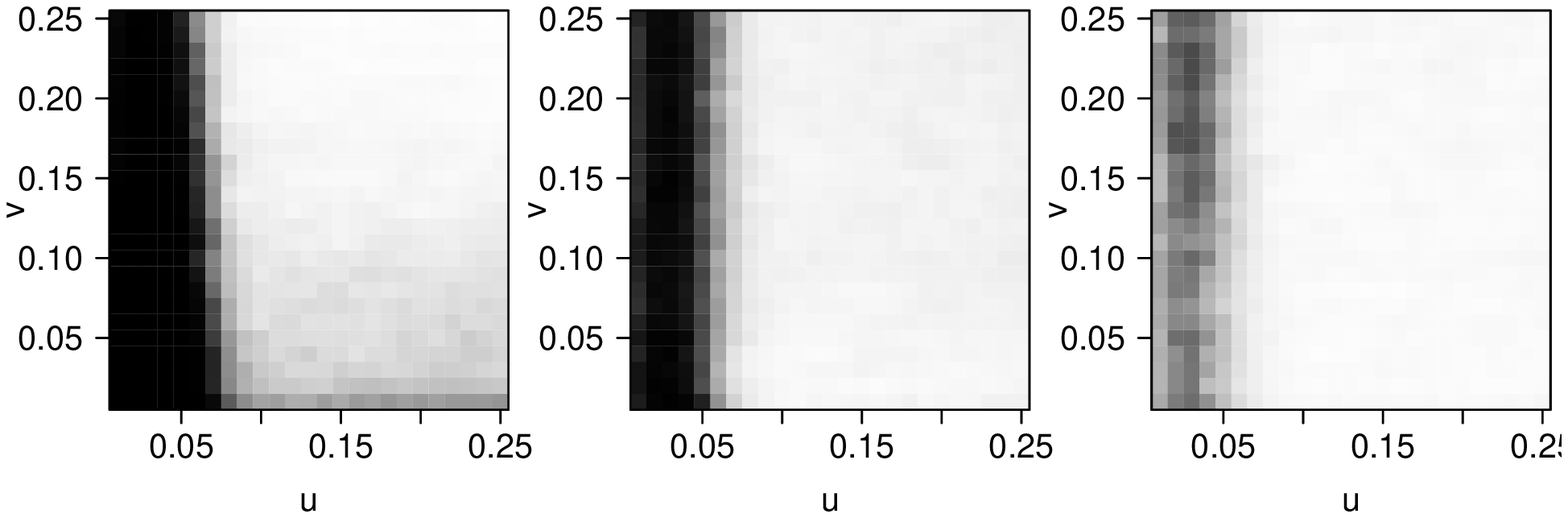}

\vspace{-1cm}
  {\tiny $PCP_3(\lambda_2(x),\beta_2=3,\sigma)$}

    \vspace{-1cm}
\includegraphics[width=7.5cm,height=4.4cm]{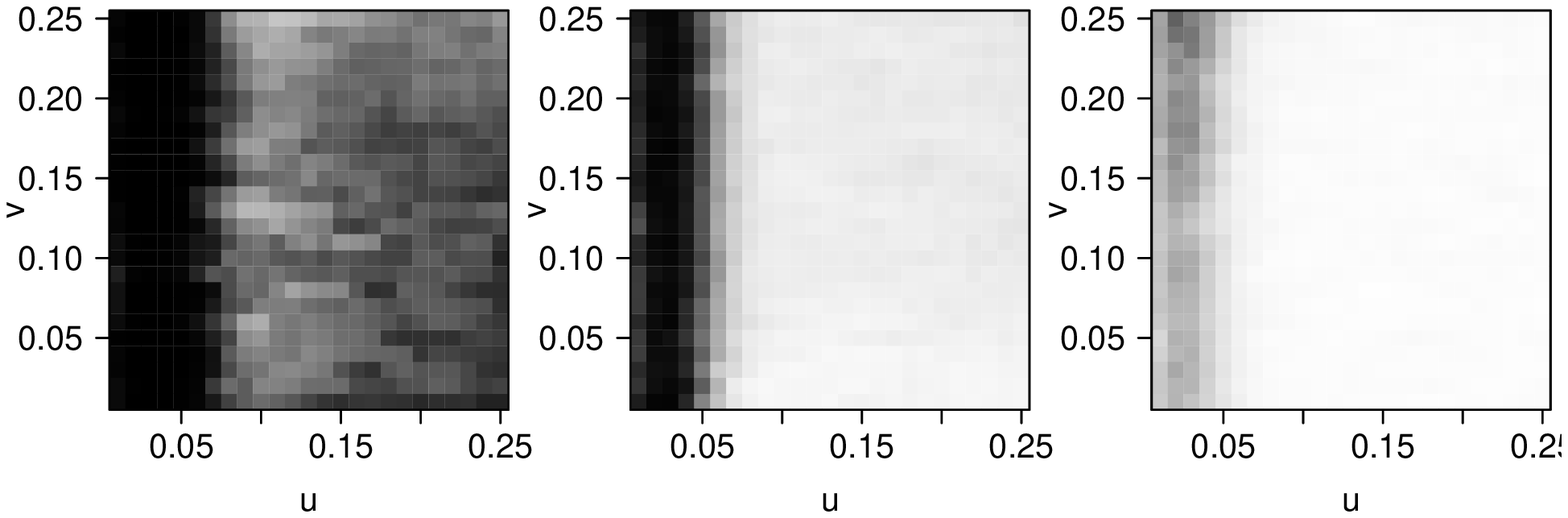}
  \end{center}

\vspace{-1.3cm}
  \caption[]{\protect\parbox[t]{5in}{Probabilities of detecting clustering using the true intensity function (first column), the non-parametric estimate with bandwidth $h=0.4$ (second column) and the parametric estimate from Model 1 (third column)
  for the Poisson cluster processes $PCP_3(\lambda(x),\beta,\sigma)$
  with $\sigma=$ 0.025.}}
  \label{fig:prob2}
\end{figure}

We analyse the detection of spatio-temporal clustering as in Section~\ref{subsec:clust}.
The probability of detecting clustering is computed for $\widehat{ g_{\lambda}}$, $\widehat{g_{\hat \lambda_p}}$ and $\widehat{g_{\hat \lambda_k}}$ from realisations of the Poisson cluster processes $PCP_3(\lambda(x),\beta,\sigma)$.
One way to get around the difficulty of estimating both the first and second-order characteristics non-parametrically consists in assuming that first-order effects operate at larger spatial scale than the second-order effects \citep{baddeley2000,diggle2007}. Thus, $\widehat{g_{\hat \lambda_k}}$ is also evaluated for fixed spatial bandwidths, $h=0.2$ and $h=0.4$.
For the processes $PCP_3(\lambda_1(x),\beta_1=1,\sigma)$ and  $PCP_3(\lambda_2(x),\beta_2=3,\sigma)$, Figure~\ref{fig:prob2} shows the probabilities of detecting clustering when using the true intensity (left), the kernel estimation with bandwidth $h=0.4$ (middle) and the parametric model 1 (right). It appears that clustering is well detected for $\widehat{g_{\lambda}}$ and $\widehat{g_{\hat \lambda_k}}$. We get a lower detection rate for $\widehat{g_{\hat \lambda_p}}$.

Finally, we calculate for each estimator the overall statistic $D$ and the related probabilities of detecting clustering. Because $D$ accumulates deviations of $\widehat g$ from the theoretical $g$ and because the bias of $\widehat{g_{\hat \lambda_k}}$ and $\widehat{g_{\hat \lambda_p}}$ may have high values, to avoid misleading results for these two estimators we rather consider the analog of $\widetilde D$ defined in Equation~(\ref{eq:dev2}).
Table~\ref{tab:power2} shows these probabilities for various scales of spatial clustering $\sigma$ and inhomogeneity.
In the simulated Poisson cluster processes, the first-order effects operate at larger spatial scale than do the second-order effects, so we can see that considering this assumption clearly improves the results and that spatio-temporal clustering can now be detected.
The estimators $\widehat{g_{\hat \lambda_k}}$ and $\widehat{g_{\hat \lambda_p}}$ are sensitive to $\sigma$  and $\beta$; their  power decrease as $\sigma$ and $\beta$ increase,
except when the bandwidth $\hat h$ is used. Indeed, it increases with $\sigma$, thus leading to less biased estimates. Conversely, the performances for a fixed bandwidth decrease with $\sigma$ because it becomes then more difficult to distinguish first- and second-order effects.
\begin{table}[h]
\centering
\caption[]{\protect\parbox[t]{5in}{Power of the integral deviation measures $\widetilde D$ according to the method used to estimate the intensity function: kernel estimation with auto-selected bandwidth $\hat h$ or fixed bandwidth $h$; parametric estimation (Model 1).}}
    \begin{tabular}{|l|c|ccc|c|}
      \hline
  Point process & $\widehat{g_{\lambda}}$ & \multicolumn{3}{|c|}{$\widehat{g_{\hat \lambda_k}}$} & $\widehat{g_{\hat \lambda_p}}$ \\
& & $\hat h$ & $h=0.2$ & $h=0.4$ & Model 1  \\
              \hline
$PCP_3(\lambda_1(x),\beta_1=1,\sigma=0.025)$ & 1.000 &0.000 &0.999 &1.000 &0.999 \\
$PCP_3(\lambda_1(x),\beta_1=1,\sigma=0.05)$  & 0.947 &0.043 &0.767 &1.000 &0.770 \\
$PCP_3(\lambda_1(x),\beta_1=1,\sigma=0.1)$   & 0.570 &0.311 &0.528 &0.997 &0.369 \\
$PCP_3(\lambda_1(x),\beta_1=2,\sigma=0.025)$ & 0.495 &0.000 &0.238 &0.600 &0.406 \\
$PCP_3(\lambda_1(x),\beta_1=2,\sigma=0.05)$  & 0.129 &0.002 &0.076 &0.242 &0.193 \\
$PCP_3(\lambda_1(x),\beta_1=2,\sigma=0.1)$   & 0.073 &0.017 &0.039 &0.097 &0.120 \\
$PCP_3(\lambda_2(x),\beta_2=3,\sigma=0.025)$ & 1.000 &0.000 &1.000 &1.000 &0.983 \\
$PCP_3(\lambda_2(x),\beta_2=3,\sigma=0.05)$  & 1.000 &0.030 &0.858 &1.000 &0.648 \\
$PCP_3(\lambda_2(x),\beta_2=3,\sigma=0.1)$   & 1.000 &0.220 &0.596 &0.994 &0.288 \\
$PCP_3(\lambda_2(x),\beta_2=5,\sigma=0.025)$ & 1.000 &0.003 &0.481 &0.801 &0.972 \\
$PCP_3(\lambda_2(x),\beta_2=5,\sigma=0.05)$  & 1.000 &0.126 &0.191 &0.334 &0.604 \\
$PCP_3(\lambda_2(x),\beta_2=5,\sigma=0.1)$   & 1.000 &0.268 &0.096 &0.107 &0.188 \\
      \hline
    \end{tabular}
\label{tab:power2}
\end{table}

\section{Discussion}
\label{sec:4}

The spatio-temporal inhomogeneous $K$-function and pair correlation function describe second-order characteristics of point processes. They can be used as an exploratory tools, e.g. for testing spatio-temporal clustering or spatio-temporal interaction \citep{diggle2010,moller2012}. Their non-parametric estimation depends on an edge correction factor and assumes that the intensity function is known.
We extended classical spatial edge correction factors to the spatio-temporal setting and
explored the influence of edge correction methods and intensity estimation on the performance of these estimators. Results on edge effects indicate that the border method performs well when the point pattern may be assumed homogeneous, otherwise, as soon as the point pattern tends to be inhomogeneous/clustered/anisotropic the translation method must be preferred.
Some preliminary tests for isotropy and Poisson assumption should be done before estimating second-order characteristics; see for instance \cite{guan2006} for testing isotropy in the spatial case or \cite{illian2008} for some discussion on the detection of anisotropy, and \cite{diggle2003} for testing the Poisson assumption. In the case of anisotropy, one then can apply adapted version of the $K$-function and pair correlation \citep{moller2012b,stoyan1994}. Note that estimators without edge correction are severely biased for values of the spatial and temporal distances larger than the ones used in this study. This is particularly so for the estimator of the STIK-function.
Results related to the intensity estimation show that the performance of the pair correlation function can be severely altered by the intensity estimate. This can be explained by over-parametrisation or over-fitting in the case of a parametric estimation of the intensity function or by the incapacity of distinguish first- and second-order effects from a single realisation of the point process in the case of a kernel-based estimation.
Non-parametric estimation of first- and second-order characteristics requires the assumption that first- and second-order effects operate at different scales.
So, an important question is what to do if first- and second-order effects operate at the same scale? It would be interesting to understand where does the confusion between the effects come from.

\section*{Appendix 1}

The proof of the unbiasedness property of the STIK-function estimator given in Equation~(\ref{eq:STIKhat1}) is based on the following Proposition, see for example \cite{moller2003}.
\begin{proposition}
Suppose that a point process $X$ defined on $W \subseteq \bR^d$ has a second-order intensity $\lambda_2$.
For a function $f : \bR^d \times \bR^d \to [0, \infty)$,
\begin{equation}
\label{eq:prop}
    \bE \lck \sum_{x,y \in X_W }^{\neq} f(x,y) \rck = \int \int f(x,y) \lambda_2(x,y) \dd x \dd y.
\end{equation}
\end{proposition}
For $W=S \times T \in \bR^2 \times \bR$ and $x=(s_x,t_x)$, $y=(s_y,t_y)$, we have from Equation~(\ref{eq:STIKhat1})
\begin{equation*}
    \widehat {K}(u,v) = \sum_{x,y \in X_W}^{\neq}
  \frac{1}{w(x,y)} \frac{1}{\lambda(x) \lambda(y)}{\bf 1}_{\lce \|s_x - s_y\| \leq u
  \ ; \ |t_x - t_y| \leq v \rce}.
\end{equation*}
Thus, for $h=x-y=(h_s,h_t)$ and $f(x,y) =  \dfrac{{\bf 1}_{\lce \|s_x - s_y\|
\leq u  \ ; \ |t_x - t_y| \leq v \rce}}{w(x,y) \lambda(x) \lambda(y)}$ we have from Equation~(\ref{eq:prop})
\begin{eqnarray*}
  \bE \lck \widehat {K}(u,v) \rck &=&  \int \int \frac{1}{w(x,y)} \frac{\lambda_2(x,y)}{\lambda(x) \lambda(y)}{\bf 1}_{\lce \|s_x - s_y\| \leq u
  \ ; \ |t_x - t_y| \leq v \rce} \Un_{W}(x) \Un_W(y) \dd x \dd y \\
  & = &  \int \underbrace{\int \frac{1}{w(x,x+h)} \Un_{W}(x) \Un_W(x+h) \dd x}_{(*)} g(h) {\bf 1}_{\lce \|h_s\| \leq u   \ ; \ |h_t| \leq v \rce}  \dd h \\
\end{eqnarray*}
For the isotropic edge correction method, denoting $w_t(x,y)=1$ if $\lck t_x-|t_x-t_y|;
\right.$ $\left. t_x+|t_x-t_y| \rck \cap T \neq \emptyset$ and 1/2 otherwise, and $|\partial b(s,r)|$ the circumference of the disc of center $s$ and radius $r$,
\begin{eqnarray*}
  (*) &=&  \dfrac{1}{|S \times T|} \int \int \frac{|\partial b(s_x,h_s)|}{|\partial b(s_x,h_s) \cap S|} \dfrac{\un_{S \times T}(s_x,t_x) \Un_{S \times T}(s_x+h_s,t_x+h_t) }{w_t(x,x+h)} \dd s_x \dd t_x \\
  & = & \frac{1}{|S \times T|} \int \frac{|\partial b(s_x,h_s)|}{|\partial b(s_x,h_s) \cap S|} \un_{S}(s_x) \Un_{S}(s_x+h_s) \dd s_x  \int \frac{\un_{T}(t_x) \Un_{T}(t_x+h_t)}{w_t(x,x+h)}
    \dd t_x \\
   & =& 1.
\end{eqnarray*}
For the modified border method,
\begin{eqnarray*}
  (*) &=& \frac{1}{|S_{\ominus u} \times T_{\ominus v}|} \int \Un_{\lce S_{\ominus u} \times T_{\ominus v} \rce}(x) \Un_W (x + h) {\bf 1}_{\lce \|h_s\| \leq u   \ ; \ |h_t| \leq v \rce} \dd x \\
  & = &   \frac{1}{|S_{\ominus u} \times T_{\ominus v}|} \int \Un_{\lce S_{\ominus u} \times T_{\ominus v} \rce}(x) \Un_{W_{\ominus h}} (x) {\bf 1}_{\lce \|h_s\| \leq u   \ ; \ |h_t| \leq v \rce} \dd x \\
  & = &   \frac{1}{|S_{\ominus u} \times T_{\ominus v}|} \int \Un_{\lce (S_{\ominus u} \cap S_{\ominus h_s}) \times (T_{\ominus v} \cap T_{\ominus h_t}) \rce}(x) {\bf 1}_{\lce \|h_s\| \leq u   \ ; \ |h_t| \leq v \rce} \dd x \hspace{3cm}\\
  & = &  \frac{1}{|S_{\ominus u} \times T_{\ominus v}|} \int \Un_{\lce S_{\ominus u} \times T_{\ominus v} \rce}(x) \dd x =  1.
  \end{eqnarray*}
For the translation method,
\begin{eqnarray*}
  (*) &=& \int \int \frac{1}{|S \cap S_{h_s}| \times |T \cap T_{h_t}|}
\Un_{S \times T}(s_x,t_x) \Un_{S \times T}(s_x+h_s,t_x+h_t) \dd s_x \dd t_x \\
& = & \frac{1}{|S \cap S_{h_s}| \times |T \cap T_{h_t}|} \int \int
\Un_{S \cap S_{-h_s}} (s_x) \Un_{T \cap T_{-h_t}}(t_x) \dd s_x \dd t_x = 1.
\end{eqnarray*}
Thus,
\begin{eqnarray*}
  \bE \lck \widehat {K}(u,v) \rck
  & = &  \int g(h) {\bf 1}_{\lce \|h_s\| \leq u   \ ; \ |h_t| \leq v \rce}  \dd h \\
  & = & \int_{-v}^v \int_{b(0,u)} g(h_s,h_t ) \dd h_s \dd h_t = K(u,v)
\end{eqnarray*}
Under the assumption of isotropy, we obtain Equation~(\ref{eq:STIKfunction1}) by using cylindrical coordinates of $g(h)$.

\clearpage

\section*{Appendix 2}

\begin{figure}[h]
 \begin{center}

  \vspace{-.6cm}
  \includegraphics[width=11.6cm,height=2.8cm]{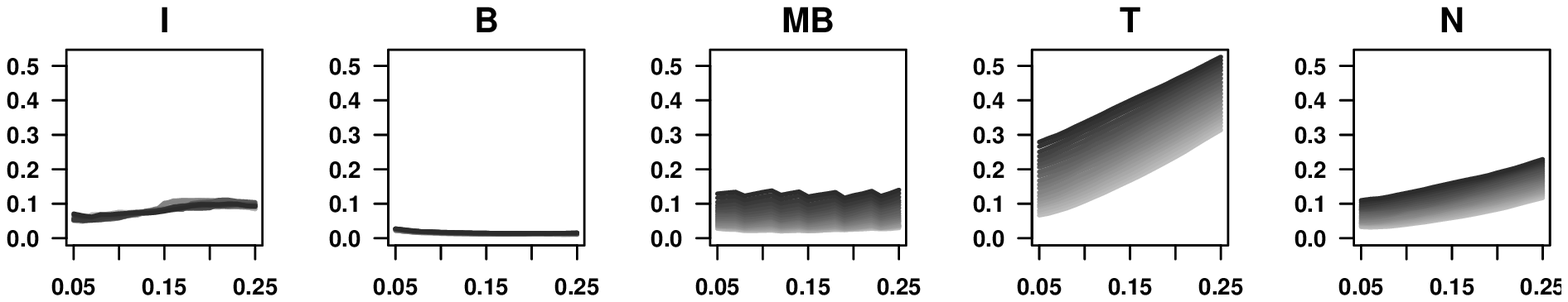}

\vspace{-.6cm}
 \includegraphics[width=11.6cm,height=2.8cm]{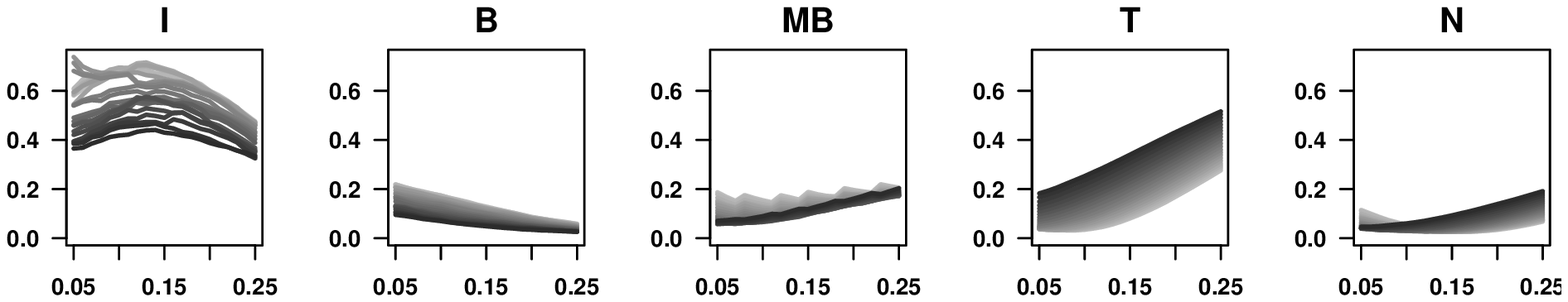}

\vspace{-.6cm}
  \includegraphics[width=11.6cm,height=2.8cm]{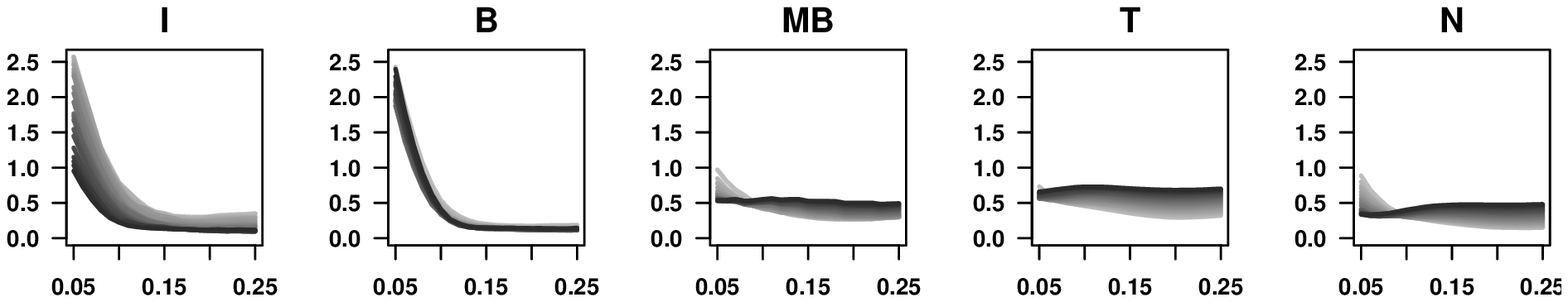}

  \vspace{-.6cm}
   \includegraphics[width=11.6cm,height=2.8cm]{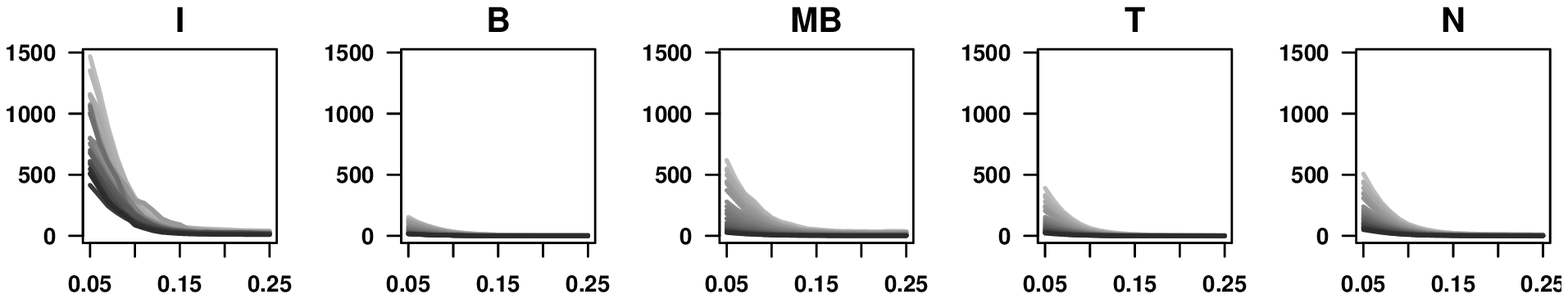}

     \vspace{-.6cm}
   \includegraphics[width=11.6cm,height=2.8cm]{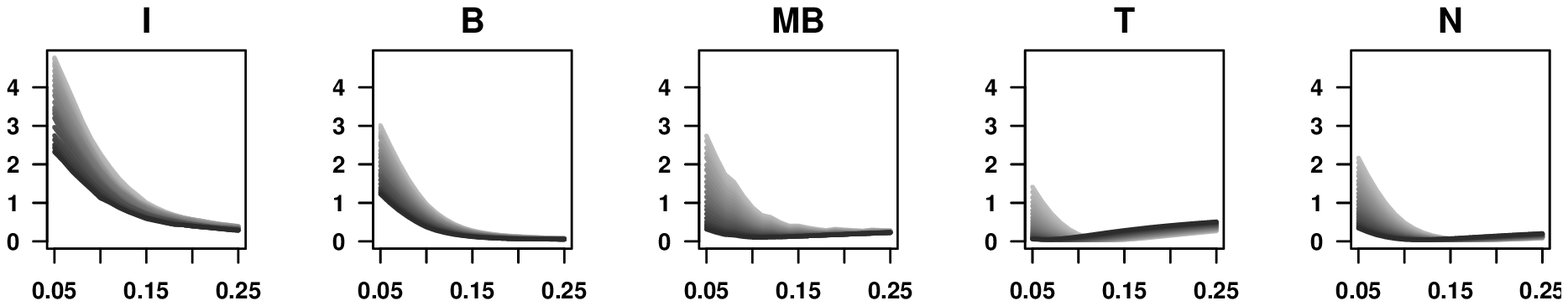}

     \vspace{-.6cm}
   \includegraphics[width=11.6cm,height=2.8cm]{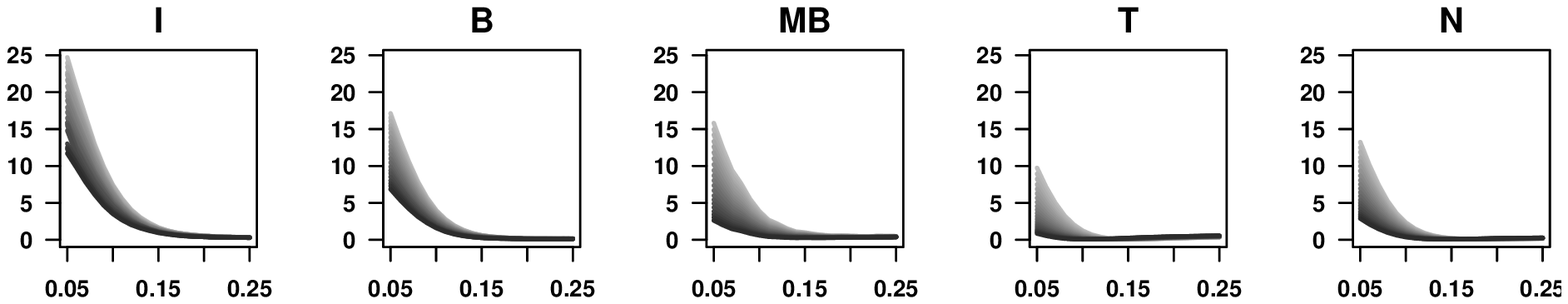}
 \end{center}

   \vspace{-.5cm}
\caption[]{\protect\parbox[t]{5in}{Mean squared error of the pair correlation function from $IPP(\lambda_1(x),\beta_1=1)$ (first row),
$PCP_1(\sigma=0.05)$ (second row), $PCP_3(\lambda_1(x),\beta_1=1,\sigma=0.05)$ (third row), $PCP_3(\lambda_2(x),\beta_2=3,\sigma=0.05)$ (fourth row),
$PCP_2(\sigma=0.2, \zeta^2=0.0625)$ (fifth row) and $PCP_2(\sigma=0.2, \zeta^2=0.5)$ (sixth row).}}
\label{fig:msePCF}
\end{figure}

Figure~\ref{fig:msePCF} shows the mean squared error of the pair correlation function estimated when using the isotropic (I), border (B), modified-border (MB), translation (T), none (N) edge correction methods. Spatial distances $u$ are in abscissa, while temporal distances $v$ are in grey: the darker, the greater, with the same range of values than $u$.

\clearpage

\begin{table}[h]
\centering
\caption[]{\protect\parbox[t]{5in}{Relative efficiency of the integral deviation measures $D$ based on the STIK-function according to the isotropic (I), border (B), modified border (MB), translation (T) edge correction factor and without edge correction (N).}}
{\scriptsize
    \begin{tabular}{|l|rrrrr|}
      \hline
        & \multicolumn{5}{|c|}{Edge correction factor} \\
             Point process & I & B & MB & T & N \\
              \hline
$HPP(\lambda)$            &  5.56 & 100.00  &     1.42  &   13.10 & 22.09 \\
$IPP(\lambda_1(x),\beta_1=1)$           &  5.30 & 100.00    &   2.13  &    7.20 & 12.73 \\
$IPP(\lambda_1(x),\beta_1=2)$           &  2.95 & 100.00    &   4.35  &    6.83 & 12.86
 \\
$IPP(\lambda_2(x),\beta_2=3)$          &  5.21  & 100.00    &  11.37   &  54.57 & 64.74\\
$IPP(\lambda_2(x),\beta_2=5)$          &  3.74 & 100.00    &   9.23   &  14.84 & 27.22 \\
$PCP_1(\sigma=0.025)$      &  0.41 &  1.15   &    0.16 &   100.00  &13.92 \\
$PCP_1(\sigma=0.05)$       &  0.20 &  0.67   &    0.15 &   100.00  & 9.26 \\
$PCP_1(\sigma=0.1)$        &  0.08 &  1.13   &    0.26 &   100.00  & 9.62 \\
$PCP_1(\sigma=0.15)$       &  0.18 &  1.63   &    0.34 &   100.00  &25.01 \\
$PCP_1(\sigma=0.2)$        &  0.30 &  4.85   &    0.99 &   100.00  &57.11 \\
$PCP_3(\lambda_1(x),\beta_1=1,\sigma=0.025)$   &  2.02 & 20.46  &    11.94 &   100.00 & 37.28 \\
$PCP_3(\lambda_1(x),\beta_1=1,\sigma=0.05)$  &  1.88  &21.61    &  15.47   & 100.00 & 30.14
 \\
$PCP_3(\lambda_1(x),\beta_1=1,\sigma=0.1)$     &  1.37 &  8.60  &     6.28  &  100.00 & 21.11 \\
$PCP_3(\lambda_1(x),\beta_1=1,\sigma=0.15)$    &  1.69 & 19.57   &   22.76  &  100.00 & 26.05 \\
$PCP_3(\lambda_1(x),\beta_1=1,\sigma=0.2)$     &  1.55  & 26.00   &   18.06  &  100.00 & 26.71 \\
$PCP_3(\lambda_1(x),\beta_1=2,\sigma=0.025)$  & 3.36 & 71.25  &   100.00  &   61.53  & 20.90 \\
$PCP_3(\lambda_1(x),\beta_1=2,\sigma=0.05)$   & 1.08 & 60.08  &   100.00  &   19.78   & 7.13 \\
$PCP_3(\lambda_1(x),\beta_1=2,\sigma=0.1)$    & 0.44 & 65.55  &   100.00   &   7.36   & 2.54 \\
$PCP_3(\lambda_1(x),\beta_1=2,\sigma=0.15)$   &  1.09 & 64.23  &   100.00  &   31.25   & 8.06 \\
$PCP_3(\lambda_1(x),\beta_1=2,\sigma=0.2)$    &  0.67 & 80.70  &   100.00  &   10.44   & 3.18 \\
$PCP_3(\lambda_2(x),\beta_2=3,\sigma=0.025)$ &  0.77& 100.00   &    1.29   &  18.44   & 8.93 \\
$PCP_3(\lambda_2(x),\beta_2=3,\sigma=0.05)$  &  0.52 &100.00   &    0.75   &  14.38   & 6.17 \\
$PCP_3(\lambda_2(x),\beta_2=3,\sigma=0.1)$   &  0.35 &100.00   &    0.65   &  18.05   & 6.78 \\
$PCP_3(\lambda_2(x),\beta_2=3,\sigma=0.15)$  &  0.25 &100.00   &    0.42   &  12.85   & 4.84 \\
$PCP_3(\lambda_2(x),\beta_2=3,\sigma=0.2)$   &  0.24 &100.00    &   0.44   &  11.02   & 4.13 \\
$PCP_3(\lambda_2(x),\beta_2=5,\sigma=0.025)$ &  0.64 &100.00    &   1.75   &   5.73   & 3.09 \\
$PCP_3(\lambda_2(x),\beta_2=5,\sigma=0.05)$  &  0.40 &100.00    &   1.96   &   5.79   & 2.70 \\
$PCP_3(\lambda_2(x),\beta_2=5,\sigma=0.1)$   &  0.17 &100.00    &   0.79   &   4.13   & 1.76 \\
$PCP_3(\lambda_2(x),\beta_2=5,\sigma=0.15)$  &  0.14 &100.00    &   0.69   &   3.64   & 1.52 \\
$PCP_3(\lambda_2(x),\beta_2=5,\sigma=0.2)$   &  0.10 &100.00    &   0.93   &   2.38   & 0.99 \\
      \hline
    \end{tabular} }
          \label{tab:devedgeK}
\end{table}

\clearpage

\begin{table}[h]
\centering
\caption[]{\protect\parbox[t]{5in}{Relative efficiency of the integral deviation measures $D$ based on the pair correlation function according to the isotropic (I), border (B), modified border (MB), translation (T) edge correction factor and without edge correction (N).}}
{\scriptsize
    \begin{tabular}{|l|rrrrr|}
      \hline
        & \multicolumn{5}{|c|}{Edge correction factor} \\
             Point process & I & B & MB & T & N \\
              \hline
$HPP(\lambda)$            &  0.79 &100.00   &   10.05  &   15.60 &23.21 \\
$IPP(\lambda_1(x),\beta_1=1)$           &  4.22 &100.00   &   24.88  &   18.11 &31.56 \\
$IPP(\lambda_1(x),\beta_1=2)$           &  5.68 & 47.31   &   60.84  &  100.00 &48.77 \\
$IPP(\lambda_2(x),\beta_2=3)$          &  0.01 & 26.07   &   23.90  &  100.00 &59.21 \\
$IPP(\lambda_2(x),\beta_2=5)$          &  0.03 & 44.02   &   41.72  &  100.00 &86.46 \\
$PCP_1(\sigma=0.025)$      &  3.80 &  9.09   &   17.97  &  100.00 &28.73 \\
$PCP_1(\sigma=0.05)$       &  1.75 &  7.71   &   12.54  &  100.00 &27.70 \\
$PCP_1(\sigma=0.1)$        &  0.32 &  4.66   &    6.49  &  100.00 &22.43 \\
$PCP_1(\sigma=0.15)$       &  0.36 &  4.96   &    3.22  &  100.00 &33.64 \\
$PCP_1(\sigma=0.2)$        &  0.39 &  6.05   &    5.96  &  100.00 &42.57 \\
$PCP_2(\sigma=0.025, \zeta^2=0.0625)$  & 17.26 & 28.64   &   44.93  &  100.00 &55.38 \\
$PCP_2(\sigma=0.05, \zeta^2=0.0625)$   & 10.76 & 24.72   &   43.67  &  100.00 &51.26 \\
$PCP_2(\sigma=0.1, \zeta^2=0.0625)$    &  5.42 & 18.40   &   24.01  &  100.00 &42.46 \\
$PCP_2(\sigma=0.15, \zeta^2=0.0625)$   &  2.59 & 12.72   &   19.71  &  100.00 &36.33 \\
$PCP_2(\sigma=0.2, \zeta^2=0.0625)$    &  1.56 & 10.67   &   20.48  &  100.00 &32.35 \\
$PCP_2(\sigma=0.025, \zeta^2=0.25)$  & 23.34 & 32.08   &   50.91  &  100.00 &58.44 \\
$PCP_2(\sigma=0.05, \zeta^2=0.25)$   & 14.71 & 27.28   &   44.38  &  100.00 &53.93 \\
$PCP_2(\sigma=0.1, \zeta^2=0.25)$    &  6.21 & 20.36   &   27.30  &  100.00 &45.17 \\
$PCP_2(\sigma=0.15, \zeta^2=0.25)$   &  4.24 & 14.07   &   18.11  &  100.00 &39.83 \\
$PCP_2(\sigma=0.2, \zeta^2=0.25)$    &  2.15 & 12.97   &   19.72  &  100.00 &35.26 \\
$PCP_2(\sigma=0.025, \zeta^2=0.5)$  & 26.39 & 34.33   &   52.15  &  100.00 &60.31 \\
$PCP_2(\sigma=0.05, \zeta^2=0.5)$   & 18.39 & 31.36   &   50.26  &  100.00 &57.38 \\
$PCP_2(\sigma=0.1, \zeta^2=0.5)$    & 13.32 & 26.29   &   35.91  &  100.00 &51.86 \\
$PCP_2(\sigma=0.15, \zeta^2=0.5)$   &  5.50 & 18.39   &   27.61  &  100.00 &46.79 \\
$PCP_2(\sigma=0.2, \zeta^2=0.5)$    &  4.44 & 16.93   &   19.81  &  100.00 &42.37 \\
$PCP_2(\sigma=0.025, \zeta^2=0.75)$  & 26.53 & 36.53   &   54.16  &  100.00 &61.38 \\
$PCP_2(\sigma=0.05, \zeta^2=0.75)$   & 23.29 & 33.84   &   54.01  &  100.00 &59.30 \\
$PCP_2(\sigma=0.1, \zeta^2=0.75)$    & 17.11 & 29.93   &   41.53  &  100.00 &54.43 \\
$PCP_2(\sigma=0.15, \zeta^2=0.75)$   &  6.87 & 21.24   &   36.98  &  100.00 &50.36 \\
$PCP_2(\sigma=0.2, \zeta^2=0.75)$    &  9.50 & 20.00   &   24.96  &  100.00 &47.72 \\
$PCP_2(\sigma=0.025, \zeta^2=1)$  & 32.54 & 36.99   &   53.98  &  100.00 &61.74 \\
$PCP_2(\sigma=0.05, \zeta^2=1)$   & 24.98 & 34.67   &   53.14  &  100.00 &59.70 \\
$PCP_2(\sigma=0.1, \zeta^2=1)$    &  3.52 & 31.16   &   44.69  &  100.00 &56.02 \\
$PCP_2(\sigma=0.15, \zeta^2=1)$   &  8.60 & 23.10   &   42.20  &  100.00 &52.47 \\
$PCP_2(\sigma=0.2, \zeta^2=1)$    & 10.22 & 22.72   &   30.69  &  100.00 &49.94 \\
$PCP_3(\lambda_1(x),\beta_1=1,\sigma=0.025)$   &  6.70 &  7.29   &   74.58  &  100.00 &46.59 \\
$PCP_3(\lambda_1(x),\beta_1=1,\sigma=0.05)$    &  3.36 &  6.03   &   47.78  &  100.00 &33.69 \\
$PCP_3(\lambda_1(x),\beta_1=1,\sigma=0.1)$     &  2.71 &  5.79   &   15.02  &  100.00 &29.54 \\
$PCP_3(\lambda_1(x),\beta_1=1,\sigma=0.15)$    &  0.36 &  6.80   &   19.07  &  100.00 &30.72 \\
$PCP_3(\lambda_1(x),\beta_1=1,\sigma=0.2)$     &  1.99 &  7.53   &   23.50  &  100.00 &26.82 \\
$PCP_3(\lambda_1(x),\beta_1=2,\sigma=0.025)$  & 12.90 & 16.90   &   59.70  &  100.00 &45.67 \\
$PCP_3(\lambda_1(x),\beta_1=2,\sigma=0.05)$   & 12.63 & 41.33   &   76.29  &  100.00 &47.71 \\
$PCP_3(\lambda_1(x),\beta_1=2,\sigma=0.1)$    & 11.78 & 97.97   &   69.36  &  100.00 &43.24 \\
$PCP_3(\lambda_1(x),\beta_1=2,\sigma=0.15)$   &  5.94 & 47.25   &   75.73  &  100.00 &38.38 \\
$PCP_3(\lambda_1(x),\beta_1=2,\sigma=0.2)$    &  4.76 &100.00   &   31.79  &   62.56 &21.84 \\
$PCP_3(\lambda_2(x),\beta_2=3,\sigma=0.025)$ &  0.24 &100.00   &    4.82  &    9.05 & 5.62 \\
$PCP_3(\lambda_2(x),\beta_2=3,\sigma=0.05)$  &  0.28 &100.00   &    2.39  &    6.36 & 3.67 \\
$PCP_3(\lambda_2(x),\beta_2=3,\sigma=0.1)$   &  0.06 &100.00   &    1.67  &    6.87 & 3.43 \\
$PCP_3(\lambda_2(x),\beta_2=3,\sigma=0.15)$  &  0.05 &100.00   &    0.90  &    3.88 & 1.93 \\
$PCP_3(\lambda_2(x),\beta_2=3,\sigma=0.2)$   &  0.06 &100.00   &    0.74  &    3.13 & 1.56 \\
$PCP_3(\lambda_2(x),\beta_2=5,\sigma=0.025)$ &  0.29 &100.00   &    2.56  &    4.48 & 3.24 \\
$PCP_3(\lambda_2(x),\beta_2=5,\sigma=0.05)$  &  0.23 &100.00   &    2.13  &    3.85 & 2.54 \\
$PCP_3(\lambda_2(x),\beta_2=5,\sigma=0.1)$   &  0.10 &100.00   &    1.19  &    3.24 & 1.94 \\
$PCP_3(\lambda_2(x),\beta_2=5,\sigma=0.15)$  &  0.04 &100.00   &    0.84  &    2.50 & 1.46 \\
$PCP_3(\lambda_2(x),\beta_2=5,\sigma=0.2)$   &  0.10 &100.00   &    0.87  &    2.18 & 1.27 \\
      \hline
    \end{tabular} }
          \label{tab:devedge}
\end{table}

\clearpage

\bibliographystyle{spbasic}      
\bibliography{Gabriel_MCAP}   

\end{document}